 \documentclass[]{IEEEtran}
\usepackage{xcolor}
\usepackage{soul}
\usepackage{booktabs,siunitx}
\usepackage[cmex10]{amsmath}
\usepackage{amsfonts}
\usepackage{graphicx,caption,subcaption} 
\usepackage{sidecap}
\usepackage{mathtools}
\usepackage{enumerate}

\usepackage{amsmath}%
\usepackage{MnSymbol}%
\usepackage{wasysym}%
\newtheorem{theorem}{Theorem} 
\newtheorem{lemma}{Lemma}

\usepackage{soul}

\definecolor{lightgray}{gray}{0.9}

\begin{document}
\title{Output Feedback Controllers Based on  a Bank of  High-Gain Observers: Robustness Analysis Against Measurement Noise}

\author{Kasra~Esfandiari
	and~ Mehran~Shakarami
	\thanks{K. Esfandiari is   with the Center for Systems Science, Yale University, New Haven, CT,  USA e-mail: (kasra.esfandiari@yale.edu).}
	\thanks{M. Shakarami  {is with the Engineering and Technology Institute Groningen,
			University of Groningen, 9747 AG Groningen, The Netherlands (e-mail:
			m.shakarami@rug.nl).}} 
}

\markboth{ }%
{Esfandiari \MakeLowercase{\textit{et al.}}: Bare Demo of IEEEtran.cls for IEEE Journals}

\maketitle

\begin{abstract} This paper analyzes   output feedback control of a class of unknown nonlinear systems in the presence of measurement noise using multiple high-gain observers (MHGO). It is well-known that single high-gain observers (HGO) are not able to provide satisfactory performance when the system output is contaminated by noise. More specifically, there is a trade-off between the   convergence speed of state estimation and the bound of steady estimation error in HGO when the output measurement is contaminated by noise. In the presented scheme,  the output feedback controller utilizes the  state estimation obtained from an appropriate combination of information provided by a bank of HGOs. The proposed strategy is capable of  mitigating  the destructive effects of measurement noise and speeding up the convergence process, and it does that because it introduces an extra   design parameter.  
	The   performance recovery capabilities of MHGO-based controllers and the stability of the closed-loop system  are     discussed. 
	 Simulations are performed  on an underwater vehicle system and a mechanical  system to evaluate the performance of the MHGO-based controller. Furthermore,  a detailed  comparison between the MHGO-based controller and controllers based on  conventional HGO, HGO with switching gain, and multi-observer {approach} is provided, which shows   the superiority of the MHG-based controller over the other methods. 
	\end{abstract}

\section{Introduction}

 \IEEEPARstart{S}{tates}
 of systems have a prominent role in control theory,
and many different strategies are developed  by using them.
Since all of the system states are not measurable in
practice, different observation schemes are presented in the control literature. 
However,  most of the previous studies   were confined  to the systems with noise-free output to simplify the understudied problem. This assumption is not realistic since    measurements are mostly  contaminated by  noise;   not only does this cause  unsatisfactory performance, but  it may also push the closed-loop system into instability.
Hence, it is necessary to investigate the robustness   of   observers and observer-based controllers \cite{RobustnessObs}. If   a priori knowledge exists about the plant, Kalman filter  is known as a   powerful tool for estimation purposes \cite{Kalman}, \cite{Cybern2018}.  
In \cite{Cybern2017}, a  fusion estimation algorithm  is presented in terms of linear matrix inequalities for a class of uncertain linear systems.  This scheme is based on the assumptions that   the uncertain part of the plant    satisfies  certain conditions and multiple sensors   measure  the output. 
In   \cite{obs1}, an adaptive  observer  is designed for  a class of nonlinear systems  with known dynamics and noisy measurements, and the relation between the observation error and bound of measurement noise is derived. 
However, the assumption of availability of a priori  knowledge about system dynamics is not always valid. Moreover,  the  control problem, which is more challenging than the state estimation problem,  has   remained intact in the works above. 


 

On the other hand, high gain observers (HGOs) are well-known as powerful structures for state estimation  of nonlinear systems. These observers are capable of handling system uncertainties and providing fast and accurate estimations if their gains are chosen sufficiently large \cite{Khalil}. For control purposes, it has been shown that by feeding sufficiently fast HGO-based   state estimations into a globally bounded controller, the output feedback controller can recover the performance of the state feedback controller \cite{atassi1999separation}.  HGOs, having  these nice features,  have attracted a great deal of attention in the past few decades and have been widely used in systems and control theory \cite{cyb,4,    Mehran}. However, conventional HGOs with large gains  yield state estimations with   substantial  over/undershoots  in the transient response, known as the peaking phenomenon.  Such  behavior may result in a closed-loop system with a finite escape time, and in turn, might destabilize  the overall plant \cite{esfandiari1992output}.  {In some works,  intelligent strategies, (e.g., fuzzy systems, neural networks, etc.) are employed to estimate the system uncertainties and    that  approximation is fed into the dynamical equation of the HGO-based structure} \cite{fuzzyhgo1}, \cite{fuzzyhgo2}, \cite{ICEE2015}.  {     Although such structures may be   applicable to a wide class of systems, they do not necessarily provide a nice transient response, specially when the initial conditions are chosen arbitrarily. Because, it takes a relatively long time for the intelligent part to learn the system dynamics, and in turn, these approaches result in an oscillatory response which is an inherent drawback of single adaptive models/structures} \cite{Narendra},  \cite{iet},  \cite{NN}. 
{ On the other hand, in the past few decades, it has been shown that   multiple  model-based typologies  are capable of providing  parameter/state estimations with  improved transient response }  \cite{NaderCDC}, \cite{Switch},  \cite{postoyan2015multi}, \cite{mulobs},  \cite{CDC}, \cite{second}. { In these approaches, multiple models are run simultaneously, and the final estimation is obtained either by switching between different models} \cite{Switch},   \cite{postoyan2015multi}, \cite{mulobs} {or by   combining the available information} \cite{CDC, second}.


 
In addition to the peaking phenomenon, another problem with the conventional HGOs is sensitivity to measurement noise \cite{Khalil}. That is because the basic idea behind the conventional   HGOs is to differentiate the system output   to get   estimations of immeasurable states. Thus, the performance  of   HGO-based structures   should be evaluated with   extra attention since the effects of measurement noise will be   {greatly amplified } by differentiating the system output.  The impact   of measurement noise on the state estimations of HGOs is discussed in \cite{khalil1}, and it has been shown that the gain of observer should not be selected too   large or too small. In general,  there is a trade-off between measurement noise sensitivity and the   convergence rate of  state estimation \cite{tradeoff}. 
However, in HGO-based feedback controllers,   sufficiently fast reconstruction of system states is a must  (before that the system states leave the region of attraction) \cite{esfandiari1992output}.  
In \cite{khalilsw2009}, a new  HGO structure is proposed for nonlinear systems in the presence of noise. In this approach,   a large gain is employed, initially,  to estimate system states fast; then the observer gain is switched to a smaller value to get a better steady state behavior.  Although the basic idea behind this observation strategy is valuable, determining  the switching time  and the   transient peaks may become challenging.

{As motivated above, we will investigate output feedback control problem of nonlinear systems in the presence of measurement noise.  In this regard,   a bank of HGOs are utilized for state estimation purposes, which enables us to improve the transient response of conventional HGOs. The employed observation strategy, MHGO,  uses  all   the information gathered from various  observers simultaneously, and a weighted summation of these observations is considered as the final estimation. 
The main contributions of the paper can be summarized as follows: }
\begin{itemize}
	\item { It is shown  that there exist some weights enabling us to estimate the system states accurately and to speed up the estimation process in the presence of noise. This re-parameterization introduces an extra design parameter to the problem; hence the need for a   large gain, which is required  in the conventional HGO,  is  mitigated.}
	\item {The output feedback control problem in the presence of measurement noise is addressed, and by using the Lyapunov's direct method, it is proven that   a semi-separation principle is valid when using the state estimations provided by the MHGO.}
	\item {The robustness analysis of the MHGO-based   controller  is   discussed, and the  conditions on the bound of measurement noise and observer gain are derived. }
	\item {Capabilities of the MHGO-based controller in recovering performance of the state feedback controller are shown, and its  supremacy of with respect to   controllers based on conventional HGO, HGO with switching gains, and multi-observer approach are provided  via simulations. }
\end{itemize}

The remainder of this paper is organized as follows: The system equation and the problem under consideration are stated in Section II.  Section III includes some preliminaries about convex sets and the HGO as well as the structure of the MHGO and   comments on its performance. The key results on robustness analysis of the closed-loop  system when the MHGO-based estimations are fed into a controller  are presented in Section IV. Section V provides simulation results, and finally Section VI summarizes the paper. 

\section{System Description and Problem Formulation}
Consider a class of nonlinear systems in canonical form, 
\begin{equation}
\begin{split}
\dot x_1&=x_2\\
\dot x_2&=x_3\\
&\,\,\,\vdots\\
\dot x_n&=f(x,u)\\
y&=x_1+\nu(t)
\label{sys0}
\end{split}
\end{equation}
where  $x= \begin{bmatrix}
x_1&\cdots&x_n
\end{bmatrix} ^T\in\mathbb{R}^n$ represents the system state vector, and $u$ and $y\in\mathbb{R}$ denote the system input and  output, respectively.   Furthermore, $\nu(t)$ expresses the output measurement noise with an unknown upper bound of $\bar \nu$, i.e., $\|\nu(t)\|\leq\bar\nu$, and  $f(x,u)$ is an  unknown nonlinear function. To guarantee the uniqueness of the system solution,   $f(x,u)$  { is assumed to be locally Lipschitz in its arguments over the domain of interest and   zero in a compact positively invariant} set $\Sigma$ \cite{Khalil}. By defining $A$, $B$, and $C$ as 

\begin{equation*}
\begin{split}
&A=\begin{bmatrix}
0&1&0&\cdots&0\\
\vdots&\ddots&\ddots&\ddots&\vdots\\
0 &\cdots& \ddots&\ddots&1\\ 
0 &\cdots &\cdots& 0&0
\end{bmatrix}, B=\begin{bmatrix}
0\\\vdots\\0 \\ 1
\end{bmatrix},  C=\begin{bmatrix}
1 \\ 0\\\vdots\\0
\end{bmatrix}^T
\end{split}
\end{equation*}
one can rewrite the system dynamics \eqref{sys0} in the following compact form
\begin{equation}
\begin{split}
\dot{x}&=Ax+Bf(x,u)\\
y&=Cx+\nu(t)
\end{split}
\label{sys}
\end{equation}   
The system dynamics \eqref{sys} cover a wide range of practical systems including electrical systems, mechanical systems, chemical processes, etc. Moreover, many other   systems, which are not in the canonical form, can be transformed into the above standard form by employing  appropriate transformations.

It is assumed that if all the system states are measurable,  the following state feedback controller  is capable of making the closed-loop system   uniformly asymptotically stable concerning    set $\Sigma$ \cite{Khalil}, 
\begin{equation}
\begin{split}
    u&=g(x,\theta)\\
    \dot {\theta}&=h(x,\theta)
\end{split}
\label{cont}
\end{equation}
 where $g(.)$ and $h(.)$ are {locally Lipschitz }in their arguments over the  domain of interest and  {globally bounded functions} of $x$. Furthermore, let us denote an open connected subset of the corresponding region of attraction  by $\mathcal{S}$. { Note that the considered class of control signal } \eqref{cont} {covers a wide range of control inputs.  The control input can be designed using feedback linearizion approach, sliding mode technique,  any adaptive approach (conventional or intelligent), etc. Thus,   the analysis provided in the subsequent sections are valid    regardless of the way that the controller has been designed.
 	In other words, one can design a state feedback controller  separately and then   replace the system states by the MHGO-based state estimations. }
  
     The understudy control problem is more general than the stabilization of an equilibrium point. In other words, lots of control problems (e.g., regulation, tracking, etc.) can be treated by properly defining the  set $\Sigma$.    For instance, the stabilization problem of the origin is a special case of the   problem above in which $\Sigma=\{0\}$.

 Since the assumption of availability of all system states is not always feasible in practice, the aforementioned controller cannot be applied to all real-life processes. To relax this assumption, it is required to estimate the system states suitably and feed them back to the controller. However, in this case, the stability of the closed-loop system  should be investigated carefully.  
   In the subsequent sections, the assumption of availability of all system states is removed by utilizing an observer-based controller,   and the  robustness   of the closed-loop system when the system output is contaminated by  measurement noise is   analyzed.

\section{ Observation Structure} 
This section presents the structures of conventional HGO  and MHGO. In addition, a brief comparison between these two state estimation strategies are provided to elucidate more on the advantages obtained from combining observations collected from different sources/observers.

\subsection{ High-gain Observer}
The dynamical equation of a single HGO is as follows:
\begin{equation}
\begin{split}
\dot{\hat{x}}&=A\hat{x}+Bf_o(\hat x,u) +H(y-C\hat{x} )
\label{sss}
\end{split}
\end{equation}
where 
$H=[
\kappa_1 / \epsilon \quad
\kappa_2 / \epsilon^2 \quad
\cdots \quad
\kappa_n / \epsilon^n
]^T 
$ and
  $\epsilon\in (0,1]$.  Function  $f_o(x,u)$ is a nominal model of function $f(x,u)$ which is locally Lipschitz in its arguments, globally bounded in $\hat x$, and zero in $\Sigma$. In addition,  $\kappa_i$ are chosen such that the  real parts of all roots of polynomial  $P(s)=s^n+\kappa_1s^{n-1}+\cdots+\kappa_{n-1}s+\kappa_n$ lie in the open left-half plane. Such a selection ensures that $A-HC$ is  a Hurwitz matrix.

{
It is well-known that the single HGO} \eqref{sss}   { can estimate system state vector accurately by selecting sufficiently large gains. However, the classic HGO} \eqref{sss}  {suffers from two major issues: }
\begin{enumerate} [(i)]
\item  undesirable   peaks exist in the transient response of the estimated  states,   and if they are fed into the controller, they may push the system into instability

\item { when the measurement is noisy, one cannot  choose the   gain in observer} \eqref{sss} {arbitrarily large. More clearly,  selecting a large gain for observer} \eqref{sss} {may yield a large steady state error. }
 \end{enumerate}
  In the following subsection, several HGOs with suitable  initial conditions are run, and the collected state estimations  are employed to estimate the system state vector in a manner such that the aforementioned issues are mitigated. 
 
 \subsection{Multiple High-gain Observers}
In this section, the MHGO structure and its capabilities in providing reliable state estimations are presented.  Since this structure  utilizes some   properties of convex sets, it is useful to present the following lemma.

\begin{lemma}\cite{bakelman2012convex} 
Let $\mathcal{K} $ be a convex subset of a linear space. Then, any element of the convex hull  $\mathcal{K} $        of   $\{q_1,\cdots,q_N\}$, i.e., $q\in\mathcal{K}$, can be expressed as $q=\sum_{i=1}^{N}\beta_iq_i
$ where $\beta_i\in[0,1] $ are constant terms  and $\sum_{i=1}^{N}\beta_i=1$.
\end{lemma}


 In order to provide state estimations using multiple HGOs, inspired by  \cite{CDC}, the  dynamical equation  for  MHGO   strategy is considered as follows:
\begin{equation}
\begin{split}
\dot{\hat{x}}_i(t)&=A\hat{x}_i(t) +H(y(t)-C\hat{x}_i(t) )\\
\hat{x}_o(t)&=\sum\limits_{i=1}^{N}{\hat\beta}_i(t) \hat{x}_i(t)
\label{xo}
\end{split}
\end{equation} 
 where $i=1,\cdots, N$, $\hat{x}_i$ is the state estimation obtained from the $i$th observer. Besides,  ${\hat\beta}_i $ represent   estimations of constant parameters $\beta_i$, and they are calculated such that   the   equality  $\sum_{i=1}^{N}{\hat\beta}_i(t) =1$ holds. Note that  to be able to use Lemma 1, the number of observers should be larger  than the number of state variables, i.e., $N\geq n+1$.
 Regarding the parameters $\beta_i$, the following lemma is considered.
 
 \begin{lemma}
 Consider the state estimation \eqref{xo}.  Let the initial conditions $\hat{x}_i(0)$ be  chosen such that $x(0)$ lies in their convex hull. Then, there exist  some positive  constant terms $\bar f_0, \bar \nu$,  $\beta_i  $ with $\sum\nolimits_{i=1}^{N}{\beta}_i =1$ such that the state estimation error $ {e}(t) =x(t)-\sum\nolimits_{i=1}^{N}{\beta}_i  \hat{x}_i(t)$ depends on   $\epsilon \bar f_0$ and $\bar \nu /\epsilon^{n-1} $. 
 \end{lemma}

 \textit{Proof: }
 In order to prove the preceding lemma, let us use the facts that $\beta_i$ are constant terms and $\sum\nolimits_{i=1}^{N}{\beta}_i=1$,  and derive the dynamical equation of error $e(t)$ as follows:
 $$\dot{{e}} =\dot x -\sum_{i=1}^{N}{\beta}_i  \dot{\hat{x}}_i =\sum_{i=1}^{N}{\beta}_i \left(\dot x - \dot{\hat{x}}_i \right)$$ 
 By  substituting   \eqref{sys} and  \eqref{xo} into the preceding equation, one can  get
   \begin{equation}
\begin{split}
\dot{{e}}  &=(A-HC) {e} +Bf(x,u)-H\nu(t) 
\label{ei0}
\end{split}
\end{equation} 
   Now, let us define a scaled version of the estimation error as   
$  \eta  =D(\epsilon)e
$,
  where the matrix $D(\epsilon)$ is defined as follows:
  \begin{equation}
\begin{split}
 D(\epsilon)=\begin{bmatrix}
1&0& \cdots&0\\
0&\epsilon&\ddots&\vdots \\
\vdots &\ddots& \ddots&0 \\ 
0 &\cdots & 0&\epsilon^{n-1}
\end{bmatrix}
\end{split}
\label{D}
\end{equation}
Taking the time derivative of the scaled error $\eta$   and using  \eqref{ei0}, one can get
   \begin{equation}
  \begin{split}
     \dot \eta&=\frac{1}{\epsilon}A_o\eta+\frac{1}{\epsilon}H_o\nu(t) +\epsilon^{n-1}Bf(x,u)
     \end{split} 
  \label{etaodot0}
  \end{equation}
{where } $ 
H_o=-\epsilon DH=\begin{bmatrix}
-\kappa_1&
-\kappa_2& 
\cdots & 
-\kappa_n  
\end{bmatrix}^T$  {and}  
\begin{equation*}
\begin{split}
&A_o=\epsilon D\left(A-HC\right)D^{-1}=\begin{bmatrix}
-\kappa_1&1&0&\cdots&0\\
-\kappa_2&0&\ddots&\ddots&\vdots\\
\vdots &\vdots& \ddots&\ddots&1\\ 
-\kappa_n &0 &\cdots& 0&0
\end{bmatrix}
\end{split}
\end{equation*}

Consider  the Lyapunov function candidate $V_0(\eta)=\eta^TP_0\eta$ for system \eqref{etaodot0}, { where $P_0=P_0^T$ is a positive definite matrix whose the largest and the  smallest eigenvalues are denoted by  $\lambda_{\max}$ and $\lambda_{\min}$, respectively. It is assumed that $P_0$ satisfies the Lyapunov's equation, i.e.,}
\begin{equation}
    A_o^TP_0+P_0A_o=-I
    \label{lyap}
\end{equation}
 Taking the time derivative of $V_0(\eta)$ and utilizing  \eqref{etaodot0} and \eqref{lyap}, yield

    \begin{equation}
  \begin{split}
     \dot V_0(\eta)&=-\frac{1}{\epsilon} \|\eta\|^2
    +\frac{2}{\epsilon}\eta^TP_0H_o \nu(t)+2\epsilon^{n-1}\eta^TP_0Bf(x,u)
     \end{split} 
     \label{v100}
  \end{equation}
Due to the globally boundedness of function $f(x,u)$ in the domain of interest, one has $\|f(x,u)\|\leq \bar f_0$. Using this fact    and performing some basic mathematical manipulations on \eqref{v100}, one can get
    \begin{equation*}
  \begin{split}
     \dot V_0(\eta)&\leq-\frac{1}{\epsilon} \|\eta\|^2+ \left(2\epsilon^{n-1}\|P_0\|\bar f_0+\frac{2}{\epsilon} \| P_0H_o\| \bar\nu \right)\|\eta\|
     \end{split} 
  \end{equation*}
  The preceding  equation can be rewritten as 
      \begin{equation*}
  \begin{split}
     \dot V_0(\eta)&\leq-\frac{1}{2\epsilon} \|\eta\|^2-\frac{\|\eta\|}{2\epsilon}\left(\|\eta\|- \left(4\epsilon^{n}\|P_0\|\bar f_0+{4}  \| P_0H_o\| \bar\nu\right)\right) 
     \end{split} 
  \end{equation*}
This equation implies that $  \dot V_0(\eta)\leq-\frac{1}{2\epsilon} \|\eta\|^2$   as long as 
  \begin{equation*}
  \|\eta\|\geq  4\epsilon^n\|P_0\|\bar f_0  +4\|P_0H_o\| \bar\nu
  \end{equation*}
 As a result, { the set} $$\mathcal{S}_0=\{V_0(\eta)\leq \left(4\epsilon^n\|P_0\|\bar f_0 +4\|P_0H_o\| \bar\nu\right)^2 \lambda_{\max}\}$$
 is an invariant set for the {system}. 
 
 Using the fact that the initial conditions of observers,  $\hat x_i(0) $,  are selected such that $x(0) $ lies in their convex hull, there exist  constant terms $\beta_i$ such that  $x(0)=\sum\nolimits_{i=1}^{N}\beta_i\hat x_i(0)$ or equivalently $e(0)=0$ (see Lemma 1), and in turn, $\eta(0)=0$. 
 {Therefore,  the estimation error is initiated from inside of the invariant set $\mathcal{S}_0$; furthermore, we have $\lambda_{\min}\|\eta\|^2\leq V_0(\eta)\leq\lambda_{\max}\|\eta\|^2$. According to these facts, it is valid to say that $\|\eta\|\leq \sqrt{\frac{\lambda_{\max}}{\lambda_{\min}}}\left(4\epsilon^n\|P_0\|\bar f_0 +4\|P_0H_o\| \bar\nu\right)  $.   By using the preceding inequality and  
   $\|e \|=\|D^{-1}\eta \|\leq \frac{1}{\epsilon^{n-1}}\|\eta\|$,   one can   show that }
\begin{equation}
     \|e\| \leq \sqrt{\frac{\lambda_{\max}}{\lambda_{\min}}} \frac{4\epsilon^n\|P_0\|\bar f_0+4\|P_0H_o\|\bar\nu}{ \epsilon^{n-1}} 
     \label{upper0}
\end{equation} 
Consequently, there exists a bounded term $\delta(t,\epsilon \bar f_0,\bar \nu /\epsilon^{n-1})$ such that  $e(t)=\delta(t,\epsilon \bar f_0,\bar \nu /\epsilon^{n-1}) $.
It is worth mentioning that when there is no measurement noise, $\bar \nu=0$,  the ultimate estimation error bound  can become arbitrarily small by choosing small enough $\epsilon$. 
   $ \hfill \blacksquare  $
   
  It is well-known that  the stability and performance of   observer-based control strategies greatly depend  on the convergence rate of the observer. More clearly, as   will be shown later,    if the state estimation error enters an  invariant set   fast enough, the closed loop system is stable. 
  {As it was shown in Lemma 2, we have $e(t)=x(t)-\sum\nolimits_{i=1}^{N}{\beta_i  \hat{x}_i(t)}=\delta(t,\epsilon \bar f_0,\bar \nu /\epsilon^{n-1})$; hence the following equality holds.  }
  	\begin{equation}
x(t)=\sum\nolimits_{i=1}^{N}{\beta_i  \hat{x}_i(t)}+\delta(t,\epsilon \bar f_0,\bar \nu /\epsilon^{n-1})
\label{xx}
  	\end{equation} 
  	{Now, let us define $e_o=x-\hat x_o$ and substitute } \eqref{xx} and \eqref{xo} { into it. Thus, we get }

  	 $$e_o(t)=\sum\nolimits_{i=1}^{N}{\beta_i  \hat{x}_i}-\sum\nolimits_{i=1}^{N}{\hat\beta_i  \hat{x}_i}+\delta(t,\epsilon \bar f_0,\bar \nu /\epsilon^{n-1})$$
  	{By adding   $x=\sum\nolimits_{i=1}^N{ \hat\beta}_ix$ to and subtracting  $x=\sum\nolimits_{i=1}^N{   \beta}_ix$ from the right-hand side of the preceding equality} ($\sum\nolimits_{i=1}^N{ \hat\beta}_i =\sum\nolimits_{i=1}^N{ \beta}_i=1$), {  it is valid to conclude that }
 \begin{equation}
     e_o(t)=\sum\nolimits_{i=1}^{N}{\tilde \beta}_i  e_i+\delta(t,\epsilon \bar f_0,\bar \nu /\epsilon^{n-1})
     \label{eoo}
 \end{equation}
  where $\tilde \beta=\hat\beta_i-\beta_i$ {and $e_i=x-\hat x_i$}. 
    It is clear that since the final estimation error $e_o(t)$ is  the  multiplication of two estimation  errors $\tilde \beta_i$ and $e_i$, this observation error is capable of  entering the invariant set very fast. In other words, this type of problem  re-parameterization  (converting the state estimation problem into estimation of constant parameters $\beta_i$) expedites the   convergence process.  
     To obtain    estimations of   $\beta_i$, the following RLS algorithm is employed,
\begin{equation}
\begin{aligned}
\dot{\hat{\bar\beta}}&=-PE^TC^T(\tilde y_N +CE\hat{\bar\beta}),& \hat{\bar\beta}(0)&=\hat{\bar\beta}_0\\
\dot{P}&=-PE^TC^TCEP,& P(0)&=\gamma I
\end{aligned}
\label{RLS}
\end{equation}
 where $\hat{\bar\beta}=\begin{bmatrix}
\hat\beta_1 &\hat \beta_2 & \cdots &\hat \beta_
{N-1}
\end{bmatrix}^T $, {$\hat\beta_{N}=1-\sum_{i=1}^{N-1}\hat\beta_i$},    $ \tilde y_N =y -C\hat x_N$, $ I\in \mathbb{R}^{(N-1)\times (N-1)} $ is the identity matrix, and $ \gamma $ is a positive constant.  Furthermore, the    $ i $th column of $ E(t) $  is defined as $\hat{x}_N(t)-\hat{x}_i(t)$.

 \textit{Remark 1: }
 In conventional HGOs  the state estimation process can be performed fast enough by choosing a sufficiently small value for design parameter $\epsilon$. However,  this results in  large peaks in the transient response of the state estimation, known as peaking phenomenon, and makes  the ultimate state observation error  large. For MHGO, it was shown  that the speed of observer depends on the convergence rate of $\hat x_i$   and $\hat \beta_i$ (see \eqref{eoo}). On the other hand, it is well-known that the convergence rate of  individual observers \eqref{xo}, $\hat x_i$,   and the RLS algorithm \eqref{RLS},  $\hat \beta_i$,   depend  on   $\epsilon$ and $\gamma$, respectively. {Thus, the need for considering a very small value for $\epsilon$   can be relaxed.   To get the desired state estimation performance,  the parameter $\epsilon$ needs to be selected large for making the ultimate estimation error $\delta(\cdot)$ small as well as avoiding  the peaking,
and the parameter $\gamma$ should appropriately be chosen for improving the transient response and expediting the convergence rate. 
 In regard to the initial conditions $\hat{\beta}_i(0)$, if there is no a priori knowledge about how close  the initial condition of the $i$th observer ($\hat x_i (0)$) is to the system states, one can set the initial weights  equally, i.e., $\hat \beta_i(0)=\frac{1}{n+1}$. In the case that a prior knowledge exists, then we will give a higher initial weight to the closest observer.}

  \section{Robustness Analysis of MHGO in Feedback Control } 

 { It is well-known that performance of observer-based controllers are dictated by the utilized state estimation.}  Due to the advantages mentioned for the state estimation obtained from  MHGO,  such a estimation  is used for control purposes in this section, and the robustness and stability analyses  of the closed-loop system are fully discussed. In this case, one can feed the estimated system states $\hat x_o$ into the  control signal  \eqref{cont}, and get  the output feedback controller as  
  \begin{equation}
\begin{split}
    u&=g(\hat x_o,\theta)\\
    \dot {\theta}&=h(\hat x_o,\theta)
\end{split}
\label{cont1}
\end{equation}
In the sequel, it will be shown that the MHGO-based control signal \eqref{cont1} is capable of recovering the performance of the state feedback controller.  
 
 To analyze the performance of the closed-loop system, first   let us subtract \eqref{sys} from \eqref{xo} and get the error dynamics of each observer as follows:
  \begin{equation}
\begin{split}
\dot{{e}}_i(t)&=(A-HC) {e}_i(t)+Bf(x,u)-H\nu(t) 
\label{ei}
\end{split}
\end{equation} 
{  Now, }  define a scaled version of the estimation error    
$  \eta_i =D(\epsilon)\left({x-\hat x_i}\right)
$, 
where $D(\epsilon) $ is as presented  in \eqref{D}.  Taking the time derivative of the scaled error   $\eta_i$ and utilizing   \eqref{ei}, result in 
  \begin{equation}
\begin{split}
\epsilon \dot{{\eta}}_i(t)&=A_o {\eta}_i(t)+\epsilon^{n}Bf(x,u)+H_o\nu(t) 
\label{etai}
\end{split}
\end{equation} 
By using the fact that $ \hat \beta_N=1-\sum_{i=1}^{N-1}\hat \beta_i$ and the definition of $\hat x_o$  \eqref{xo}, one can   write the scaled   state   estimation error $\eta_o=D(\epsilon)\left(x-\hat x_o\right)$ as follows:
  \begin{equation}
  \begin{split}
      \eta_o&=D(\epsilon)\left(\sum_{i=1}^{N-1}\hat \beta_i\left(x-\hat x_i\right)+\left(1-\sum_{i=1}^{N-1}\hat \beta_i\right)(x-\hat x_N)\right)\\
      &=
  E_o\hat{\bar \beta}+\eta_N
  \end{split} 
  \label{etao}
  \end{equation}
where the $i$th column of $E_o$ is $\eta_i-\eta_N$. Besides, by employing \eqref{etai}, it is straightforward to show that 
\begin{equation}
   \epsilon \dot E_o= A_oE_o
   \label{Eo}
\end{equation} 
 
To get the dynamical equation of $\eta_o$, it is required to take the time derivative of \eqref{etao},  and employ \eqref{RLS}, \eqref{etai},  \eqref{Eo}. Thus,  one can get
  \begin{equation*}
  \begin{split}
     \dot \eta_o&=\frac{1}{\epsilon}A_o\eta_o-E_oPE^TC^TC\left(e_N+E\hat{\bar\beta}\right)\\
     &-E_oPE^TC^T\nu(t)+\epsilon^{n-1}Bf(x,u)+\frac{1}{\epsilon}H_o\nu(t)
  \end{split} 
  \end{equation*}
  By employing the facts that $D^{-1}\eta_o=e_o=e_N+E\hat{\bar\beta}$, $CD^{-1}=C$, and $E=D^{-1}E_o$, one can  rewrite the preceding equation as follows,
    \begin{equation}
  \begin{split}
     \dot\eta_o&=\frac{1}{\epsilon}A_o\eta_o-E_oPE_o^TC^TC\eta_o\\
     &+\left(-E_oPE^T_oC^T\right)+\frac{1}{\epsilon}H_o)\nu(t)+\epsilon^{n-1}Bf(x,u)
     \end{split} 
  \label{etaodot}
  \end{equation}
  
  Now, let us employ \eqref{sys} and  \eqref{etaodot} and write the system dynamics under the output feedback controller  \eqref{cont1}
    as follows:
  \begin{align}
\dot{x}&=Ax+Bf\left(x,g(x-D^{-1}\eta_o,\theta)\right) 
\label{sys1}
\\
  \epsilon\dot \eta_o&=A_o\eta_o-\epsilon E_oPE_o^TC^TC\eta_o\label{sys11}\\
     &+\left(-\epsilon E_oPE^T_oC^T+H_o\right)\nu(t)+\epsilon^{n}Bf\left(x,g(x-D^{-1}\eta_o,\theta)\right) \nonumber
\end{align}
   

  The obtained dynamical equations represent a system in the standard singularity perturbed form. 
 In order to analyze the closed-loop system behavior, one needs to consider the following  facts and lemma.

   \textit{Fact 1:} Because $P(t)$ is a positive definite matrix and $\dot P(t)\leq 0$ (see \eqref{RLS}), it is valid to conclude that $P(t)$ is bounded. 
   
   \textit{Fact 2:} Since the matrix $A_o $ is Hurwitz, the  dynamical equation \eqref{Eo} results in a bounded term $E_o(t)$.  In  other words, there exist positive constants $l_1$ and $\lambda $ such that  $
      \|E_o(t)\|=\|\exp(\frac{1}{\epsilon}A_ot)E_o(0)\|\leq l_1  \exp(-\frac{1}{\epsilon}\lambda t)
$.   

    $  $     


    \begin{lemma}
Consider  the nonlinear function $h(\epsilon,\bar\nu)=\frac{4\epsilon^n\bar f+2\left(a_1\epsilon+a_2\right)\bar\nu}{\epsilon^{n-1}}$ with  positive constants $\bar f$, $ a_1$, and $ a_2$. Then, there exist $\epsilon^*\in(0,1]$, $\epsilon_1^*<\epsilon^*$ and $\epsilon_2^*>\epsilon^*$ such that  for  every $\bar\nu\in[0,\bar\nu^*(\epsilon^*)]$ and a given  constant term $\bar h$, inequality  $h(\epsilon,\bar\nu)\leq\bar h$ holds for every $\epsilon\in [\epsilon_1^*, \epsilon_2^*]$. 
    \end{lemma} 
     

  \textit{Proof:}  
 To prove the lemma, we will first show that  $h(\epsilon,\bar\nu) $   has only one    minimum point at  $\epsilon^*$. Then, this fact will be utilized to prove the lemma. 
  
  To prove the first part, it will be first shown that $\frac{\partial h}{\partial \epsilon }=0$ has at most two roots. Then, all possible scenarios will be discussed in detail, and it will be concluded that $h(\epsilon)$  has exactly one minimum point. In this regard, the partial  derivative of $h(\epsilon)$ with respect to  $\epsilon$ can be taken as follows:
  \begin{equation}
      \frac{\partial h }{\partial \epsilon}=\frac{4\bar f \epsilon^{n}-2(n-2)a_1\bar \nu\epsilon-2(n-1)a_2\bar\nu}{\epsilon^{n}}
      \label{hder}
  \end{equation}
 To find extrema   of $h(\epsilon)$, one should set   $\frac{\partial h }{\partial \epsilon}=0$. Since $\epsilon\ne 0$, this is equivalent to setting the numerator of the preceding equation equal to zero, i.e., 
   $   h_1(\epsilon)= {4\bar f \epsilon^{n}-2(n-2)a_1\bar \nu\epsilon-2(n-1)a_2\bar\nu}=0$.
 In order to show that the number of roots of $h_1(\epsilon)=0$ is at most two,     $\frac{\partial h_1 }{\partial \epsilon}$ will be checked.  By performing some basic manipulations, one   can get   $\frac{\partial h_1 }{\partial \epsilon}=4n \bar f \epsilon^{n-1}-2(n-2)a_1\bar \nu=0$,  which    has only one solution at $ \epsilon=\left(\frac{2(n-2)a_1\bar \nu}{4n \bar f}\right)^{\frac{1}{n-1}}$. Hence, sign of $\frac{\partial h_1 }{\partial \epsilon}$   changes  one time (from a negative value to a positive value), and in turn,  it is valid to conclude that $h_1(\epsilon)=\frac{\partial h }{\partial \epsilon}= 0$   has  at most  two solutions. 
 
 Now, let us consider the following three  possible  cases: (i) $\frac{\partial h }{\partial \epsilon}=0$ has no solution (ii) $\frac{\partial h }{\partial \epsilon}=0$ has one solution (iii) $\frac{\partial h }{\partial \epsilon}=0$ has two distinct solutions. Let us check some properties of function $h(\epsilon)$ and show that   case (i)   results  in a   contradiction.    Since function $h(\epsilon) $ has a negative slope for small values of $\epsilon$ (see \eqref{hder}), $\lim_{\epsilon\rightarrow 0^+} \frac{\partial h }{\partial \epsilon}<0$,  this function is indeed decreasing at the beginning. On the other hand, we have    $\lim_{\epsilon\rightarrow 0^+}   h(\epsilon) =+\infty $ and         $\lim_{\epsilon\rightarrow +\infty}   h(\epsilon) =+\infty $. Thus, the  slope of this function should change its sign at some points, which contradicts with the assumption of having no solution for   $\frac{\partial h }{\partial \epsilon}=0$, i.e., case (i). 
  For case (ii), let us assume that $\epsilon_1$ denotes the root of    $ \frac{\partial h }{\partial \epsilon}=0$. Since  $\lim_{\epsilon\rightarrow 0^+} \frac{\partial h }{\partial \epsilon}<0 $ and both $\lim_{\epsilon\rightarrow 0^+}   h(\epsilon)$ and  $\lim_{\epsilon\rightarrow +\infty} h(\epsilon) $    tend to $+\infty$, $\epsilon_1$ is the minimum point of function $h(\epsilon)$. 
  On the other hand, we have  $\epsilon\in(0,1]$; thus,  in this case, the minimum occurs at $\epsilon^*=\min\{\epsilon_1,1\}$.
 For case (iii), let us denote the two distinct roots of $\frac{\partial h }{\partial \epsilon}=0$ by $\epsilon_2$ and $\epsilon_3$. Using a similar discussion presented  for  cases (i) and (ii), one can conclude that at least one of these distinct roots should be the minimum point of  function $h(\epsilon)$, e.g.,  $\epsilon_2$. 
  For the other root, i.e., $\epsilon_3$, since we assumed that $  \frac{\partial h }{\partial \epsilon}\big{|}_ {\epsilon=\epsilon_3} =0$, this point    can be a minimum or a maximum or an inflection point of function $h(\epsilon)$.  
 Because   no function can have two  consecutive minimum points (without having any maximum point in the between of them), $\epsilon_3$  cannot be a minimum point. It cannot be a maximum  point either since this assumption  contradicts with the fact that  $\lim_{\epsilon\rightarrow +\infty} h(\epsilon)=+\infty $. On the other hand,   $\epsilon_3$ is not an inflection point of   function $ h(\epsilon)$ since $\frac{\partial^2 h}{\partial\epsilon ^2 }=\frac{ 2(n-2)(n-1)a_1\bar \nu\epsilon+2(n-1)na_2\bar\nu}{\epsilon^{n+1}}\ne 0$ for all bounded values of $\epsilon$. In other words,  case (iii) does not occur. 
   

 So far, it was proven that    function $h(\epsilon)$ has exactly one  minimum at $\epsilon^*$. To find the largest possible (less conservative) upper bound  of  $\bar \nu$, for a given constant $\bar h$,   we need to check $h(\epsilon^*,\bar \nu)\leq \bar h$.  Toward this end, by performing some basic manipulations on  $ \frac{4\epsilon^{*n}\bar f+2\left(a_1\epsilon^*+a_2\right)\bar\nu}{\epsilon^{*n-1}}\leq \bar h$, one can conclude that the upper bound of noise, $\bar \nu$, should be less than or equal to $\bar\nu^*=\frac{ \epsilon^{*n-1}\bar h-4\epsilon^{*n}\bar f}{2\left(a_1\epsilon^*+a_2\right)}$. On the other hand, it is clear that if $\bar \nu\in [0,\bar\nu^*]$, then the equation $h(\epsilon,\nu)=\bar h$ has two solutions  at $\epsilon_1^*<\epsilon^*$ and $\epsilon_2^*>\epsilon^*$. Hence,  the inequality $h(\epsilon,\nu)\leq \bar h$ holds for    every $\epsilon\in [\epsilon_1^*,\epsilon_2^*]$.  
   $\hfill  \blacksquare  $

   The following theorem summarizes performance recovery of the singularly perturbed closed-loop system  in the presence of measurement noise (refer to  \eqref{sys1} and \eqref{sys11}).

\begin{theorem}
Let us consider  the dynamical system  \eqref{sys0} with the control input \eqref{cont}. If the system states are estimated using the  observer \eqref{xo} with the adaptive law \eqref{RLS}, then 
for any compact set ${\mathcal{S}_1} \subseteq \mathcal{S}$ (where $\mathcal{S}$ is   an open connected subset of the region of attraction)  and any compact set  $\mathcal{S}_2\subseteq \mathbb{R}^n$, there exist constants $\bar\nu^*$, $\epsilon_1^*$ and $\epsilon_2^*$ such that for every $\|\nu(t)\|\leq \bar\nu^*$, $\epsilon\in [\epsilon_1^*,\epsilon_2^*]$, the solution $(x,\hat x_o)$, starting in $\mathcal{S}_1\times\mathcal{S}_2$, is bounded for all $t$. Furthermore, the adaptive parameters $\hat {\bar\beta}(t)$ and the individual observers estimations $\hat x_i(t)$ are bounded as well. 
\end{theorem}

  \textit{Proof:} To prove the theorem,  a  positive  invariant set will be derived for the system dynamics; then this will be utilized to ensure    boundedness of all signals of the closed-loop system. Toward this end, let us  consider the Lyapunov function candidate $V_1(\eta_o)=\eta_o^TP_0\eta_o$ for system \eqref{sys11}. Taking the time derivative of this function and substituting   \eqref{sys11} and \eqref{lyap} into it, yield
    \begin{equation}
  \begin{split}
     \dot V_1(\eta_o)&=-\frac{1}{\epsilon}\eta_o^T \eta_o-2\eta_o^T   P_0 E_oPE_o^TC^TC   \eta_o \\
     &+2\eta_o^TP_0\left(-E_oPE^T_oC^T+\frac{1}{\epsilon}H_o\right)\nu(t)\\
     &+2\epsilon^{n-1}\eta_o^TP_0Bf(x,u)
     \end{split} 
     \label{v1}
  \end{equation}
 Due to the globally boundedness of controller $u$ in its arguments and locally Lipschitz property of $f(x,u)$, one has $\|f(x,u)\|\leq \bar f_0$ over a domain of interest $\mathcal{S}_c\subseteq \mathcal{S}$ ($\mathcal{S}_c$ will be defined  later). 
 By using the preceding inequality and performing some basic mathematical manipulations on \eqref{v1}, one can get
    \begin{equation*}
  \begin{split}
     \dot V_1(\eta_o)&\leq-\frac{1}{\epsilon}\|\eta_o\|^2-2 \eta_o^T P_0  E_oPE_o^TC^TC \eta_o \\
     &+\left(a_1+\frac{1}{\epsilon}a_2\right)\|\eta_o\|\bar\nu 
      +2\epsilon^{n-1}\|\eta_o\|\bar f
     \end{split} 
       \end{equation*}
  where $2\|P_0E_oPE^T_oC^T\|\leq a_1$,   $a_2=2 \| P_0H_o\|$, and $\bar f=\|P_0\|\bar f_0$.     With regard  to  Facts 1 and 2, one can conclude that constant  $a_1$ is a bounded term.

  Now let us define the  { compact set }
  \begin{equation}
  \mathcal{S}_3=\{V_1(\eta_o)\leq \left(2\left(\epsilon a_1+a_2\right) \bar\nu+4\epsilon^n\bar f \right)^2\lambda_{\max}\}
  \label{set}
  \end{equation}
  Outside of the above set, one has 
         \begin{equation*}
  \begin{split}
     \dot V_1(\eta_o)&\leq-\frac{1}{2\epsilon}\|\eta_o\|^2-2 \eta_o^T P_0  E_oPE_o^TC^TC \eta_o 
     \end{split} 
     \label{aa}
  \end{equation*}
  By using Fact 2 and performing some basic manipulations on the preceding inequality,  one can get
      \begin{equation}
  \begin{split}
     \dot V_1(\eta_o)&\leq-\frac{1}{2\epsilon}\|\eta_o\|^2+l_2\exp(-\frac{2}{\epsilon}\lambda t)\| \eta_o\|^2  \ 
     \end{split} 
     \label{v22}
  \end{equation}
  where $2l_1^2 \|P_0\| \|P\|\|C^TC\|\leq l_2$. Note that by utilizing Facts 1 and 2, it is straightforward to conclude that the positive constant $l_2$ is  bounded.
    By using inequality $\lambda_{\min}\|\eta_o\|^2\leq V_1(\eta_o)\leq \lambda_{\max}\|\eta_o\|^2$ and  \eqref{v22}, one has 
           \begin{equation}
  \begin{split}
     \dot V_1(\eta_o)&\leq\left(-\frac{1}{2\epsilon\lambda_{\max}}+\frac{l_2}{\lambda_{\min}}\exp(-\frac{2}{\epsilon}\lambda t)\right)V_1( \eta_o) \ 
     \end{split} 
     \label{v222}
  \end{equation}
  where $\lambda_{\max}$ and $\lambda_{\min}$ denote the largest and smallest eigenvalues of the matrix $P_0$. Taking integral over \eqref{v222}, { results  in }
         \begin{equation*}
  \begin{split}
      V_1(t)&\leq V_1(0)\exp\left(-\frac{t}{2\epsilon\lambda_{\max}}\right)\exp\left(\frac{l_2\epsilon}{2\lambda\lambda_{\min}}\left(1-\exp(-\frac{2}{\epsilon}\lambda t)\right)\right)
     \end{split} 
  \end{equation*}
  Since there exists a positive constant   $l_3$ such that $\exp\left(\frac{l_2\epsilon}{2\lambda\lambda_{\min}}\left(1-\exp(-\frac{2}{\epsilon}\lambda t)\right)\right)\leq l_3$, one can get
        \begin{equation}
  \begin{split}
      V_1(t)&\leq V_1(0)l_3\exp\left(-\frac{t}{2\epsilon\lambda_{\max}}\right)
     \end{split} 
     \label{v3}
  \end{equation}
 Therefore, if $\eta_o$ is outside of the compact set $\mathcal{S}_3$ \eqref{set}, there exists a finite time $T(\epsilon)$ after which  $\eta_o$  will enter that  set. To obtain a closed form for $T(\epsilon)$, the  preceding inequality can be utilized to get
\begin{equation}
T(\epsilon)={4\epsilon}\lambda_{\max} \ln{\left(\frac{\sqrt{  V_1(0)l_3}}{4\epsilon^n\bar f\sqrt{\lambda_{\max}}}\right)}
\label{T}
\end{equation} 

{On the other hand, as long as the scaled state estimation error is inside $\mathcal{S}_3$, $\|x-\hat x_o\|=\|D^{-1}\eta_o\|\leq \frac{1}{\epsilon^{n-1}}\|\eta_o\|$ satisfies the following inequality}
\begin{equation}
     \|x-\hat x_o\| \leq \sqrt{\frac{\lambda_{\max}}{\lambda_{\min}}}h(\epsilon,\bar \nu)
     \label{upper}
\end{equation} 
  where $h(\epsilon,\bar \nu)=\frac{4\epsilon^n\bar f+2\left(a_1\epsilon+a_2\right)\bar\nu}{\epsilon^{n-1}}$. Hence, we showed that $\eta_o(t)$ is bounded; however,  the provided analysis were based on the assumption that $x(t) \in \mathcal{S}_c$. 
  In the sequel, the analysis of this part is divided into two steps. First,   we will   ensure that when  $x(t) $ starts from inside of the  set $ \mathcal{S}_1\subseteq \mathcal{S}_c$,  $\eta_o(t)$ will enter the set $  \mathcal{S}_3$ before that  $x(t)$ leaves $  \mathcal{S}_c$, i.e., the provided analysis for $\eta_o(t)$ is valid during this time interval.   Second, it will be shown  that  $  \mathcal{S}_c\times  \mathcal{S}_3$ is a positive invariant set, and in turn,  $x(t)$ and $\eta_o(t)$ will remain inside the  set  $  \mathcal{S}_c\times  \mathcal{S}_3$ thereafter. 

Since the system dynamics are in the form of standard singularly perturbed systems \cite{Khalil},  let us substitute  $\eta_o=0$ into \eqref{sys1} and get     
   \begin{equation}
\begin{split}
\dot{x}&=Ax+Bf(x,g(x,\theta))
 \end{split}
 \label{state}
\end{equation}
 It is obvious that the  reduced system  is identical to the system under the state feedback controller \eqref{cont}, and in turn uniformly asymptotically stable with respect to the positively invariant set $\Sigma$. According Lyapunov's converse Theorem \cite{Khalil}, there exists a Lyapanuv's function $V_2(x)$ and positive definite functions $U_1(x)$, $U_2(x)$, and $U_3(x)$ for system \eqref{state} such that 
 \begin{equation}
 \begin{split}
 V_2(x)&=0 \iff x\in\Sigma\\
 U_1(x)&\leq V_2(x)\leq U_2(x)\\
   {\dot V_2(x)} &\leq -U_3(x)\\
 \lim_{x\xrightarrow{}\partial\mathcal{S}}U_1(x)&=\infty
 \end{split}
  \end{equation}
where $\mathcal{S}$ is an open connected subset of the  region of attraction; moreover there exists   $c\geq \max_{x\in \mathcal{S}_1} V_2(x)$ such that     $\mathcal{S}_1\subseteq\mathcal{S}_c=\{V_2(x)\leq  c\}\subseteq\mathcal{S}$.

Since the nonlinear function $f(x,u) $ is locally Lipschitz function and $u=g(x-D^{-1}\eta_o)$ is globally bounded over the set of interest  $\mathcal{S}_c$, one has
\begin{equation*}
\|\dot x\|=\|Ax+Bf(x,g(\cdot))\|\leq a_5
\end{equation*} 
where $a_5>0$ is a constant term. Taking integral over both sides of the preceding equation and using the fact that $\|\int_0^t\dot xd\tau\|\leq\int_0^t\|\dot x\|d\tau$, yields $ \| x(t)-x(0)\| \leq a_5t $. Moreover, we have $x(0)\in \mathcal{S}_1\subseteq \mathcal{S}_c$; thus the    inequality $ \| x(t)-x(0)\| \leq a_5t $ implies that there exists $T_1$ such that $x(t)$ is inside the set $\mathcal{S}_c$ as long as  $t\leq T_1$.  On the other hand,    since $T(\epsilon)$  tends to zero as $\epsilon\rightarrow 0$ (refer to \eqref{T}), there exists a constant term $\epsilon_3^*$ such that for $\epsilon\leq \epsilon_3^*$ we have $T(\epsilon)\leq T_1$.   In other words, the   scaled state estimation error $\eta_o(t)$ enters the set $ \mathcal{S}_3$ fast enough before that the system states $x(t)$ leave the set $ \mathcal{S}_c$. 

In the next step, it will be shown that if $(x,\eta_o)$ lies inside of the set $\mathcal{S}_c\times \mathcal{S}_3$, this pair will always remain there. In other words, $\mathcal{S}_c\times \mathcal{S}_3$ is  a positive invariant set. 
In this regard, by using the   Lyapunov's function $V_2(x)$ for closed-loop system \eqref{sys1}, one can get
\begin{equation*}
\begin{split}
        \dot V_2(x)&\leq -U_3(x)+\frac{\partial V_2}{\partial x}B\left(f(x,D^{-1}\eta_o)-f(x)\right)
\end{split}
\end{equation*}  
   Due to the  fact that  $\|\frac{\partial V_2}{\partial x}\|\leq a_3$ and the Lipschitz property of   function $f(\cdot)$      over the domain of interest, it is valid to conclude that 
  \begin{equation*}
\begin{split}
        \dot V_2(x)&\leq -U_3(x)+a_3a_4\|D^{-1}\eta_o \|
\end{split}
\end{equation*} 
where $a_4$ denotes the Lipschitz constant.   Inside of the set $\mathcal{S}_c\times \mathcal{S}_3$  the presented  upper bound  in \eqref{upper} is valid; thus by utilizing this upper bound and  $\|D^{-1}\eta_o\|=\|x-\hat x_o\|$, the preceding inequality  can be rewritten {as follows:}
  \begin{equation*}
\begin{split}
        \dot V_2(x)&\leq -U_3(x)+\sqrt{\frac{\lambda_{\max}}{\lambda_{\min}}}a_3a_4h(\epsilon,\bar\nu)
\end{split}
\end{equation*} 
{By setting} $\bar h=\sqrt{\frac{\lambda_{\min}}{\lambda_{\max}}}\frac{1}{a_3a_4}\min_{x\in\partial \mathcal{S}_c}U_3(x)$ in Lemma 3, we get $\dot V_2(x)\leq 0$ for $\bar\nu\in[0,\bar\nu^*_1]$ and  $\epsilon\in[\epsilon_4^*,\epsilon_5^*]$.
On the other hand,  $\epsilon_3^*$ (obtained earlier for ensuring $T(\epsilon)\leq T_1$ for all $\epsilon\leq \epsilon_3^*$) gives us an upper bound for the measurement noise, i.e., $\bar\nu_2^*=\frac{ \epsilon_3^{*n-1}\bar h-4\epsilon_3^{*n}\bar f}{2\left(a_1\epsilon_3^*+a_2\right)}$. Thus, the  parameters $\bar\nu^*$, $\epsilon_1^*$, and $\epsilon_2^*$ (used in the theorem) can be defined as  $\bar\nu^*=\min\{\bar\nu_1^*,\bar\nu_2^*\}$,  $\epsilon_1^*=\epsilon_4^*$, and $\epsilon_2^*=\min\{\epsilon_3^*,\epsilon_5^*\}$.
Note that  we showed that $  \dot V_1(\eta_o)\leq 0$ and $  \dot V_2(x)\leq 0$ for all $(x(t), \eta_o(t))\in \mathcal{S}_c\times \mathcal{S}_3$. Hence,  the set $ \mathcal{S}_c\times \mathcal{S}_3 $ is a positive invariant set.

In summary, we proved that if $x(0)\in\mathcal{S}_1\subseteq\mathcal{S}_c$ and $\eta_o$ is outside of the set $\mathcal{S}_3$, then  $x(t)$ and  $\eta_o(t)$  will enter   the  set $\mathcal{S}_c\times\mathcal{S}_3$    after $T(\epsilon) $ units of time and will remain there for  $t> T(\epsilon)$. This means that the solution $(x(t),\hat x_o(t))$ is bounded. 

To complete the proof and ensure boundedness of the other signals of the closed-loop system, it is needed to guarantee that $\hat x_i, \hat{\bar\beta}\in \mathcal{L}_\infty$. Toward this end, first  let us show that each   observer yields a bounded state estimation vector, i.e.,  $\hat x_i\in\mathcal{L}_\infty$. In this regard,  one can rewrite  dynamical equation  \eqref{xo} as
\begin{equation}
    \dot{\hat{x}}_i(t)=(A-HC)\hat{x}_i(t) +H\left(C x(t)+\nu(t)\right)
    \label{xii}
\end{equation}
It was proven earlier that  $x(t)$ belongs to $\mathcal{L}_\infty$; moreover $\nu(t)$ is    bounded  as well.    Therefore,  equation \eqref{xii} represents a linear system with Hurwitz matrix $A-HC$ and  bounded input   $C x(t)+\nu(t)$, and in turn, $\hat x_i\in\mathcal{L}_\infty$.

To prove that $\hat{\bar \beta}\in \mathcal{L}_{\infty}$,   let us take integrate over  \eqref{RLS} and get ${\hat{\bar\beta}}(t)-{\hat{\bar\beta}}(0)=-\int_{0}^{t}PE^TC^T(\tilde y_N +CE\hat{\bar\beta})d\tau$, and   in turn, one has  
$
    \|{\hat{\bar\beta}}(t)\|\leq\|{\hat{\bar\beta}}(0)\|+\|\int_{0}^{t}PE^TC^T(\tilde y_N +CE\hat{\bar\beta})d\tau\|
 $. 
Since $C\left(\tilde x_N+E\hat{\bar\beta}\right)+\nu(t)=C\tilde x_o+\nu(t)$, the preceding inequality can be rewritten as
\begin{equation}
    \|{\hat{\bar\beta}}(t)\|\leq\|{\hat{\bar\beta}}(0)\|+\int_{0}^{t}\|PE^TC^T\left(C\tilde x_o+\nu\right)\|d\tau
    \label{int}
\end{equation}
As it was shown earlier $\tilde x_o=D^{-1}\eta_o\in \mathcal{L}_\infty$; thus $ C \tilde x_o+\nu$   belongs to $\mathcal{L}_\infty$.  Using the preceding equality, Facts 1 and 2, and \eqref{int}, one has
$ \|{\hat{\bar\beta}}(t)\|\leq\|{\hat{\bar\beta}}(0)\|+a_6\int_{0}^{t}  \exp(-\frac{1}{\epsilon}\lambda \tau) d\tau$, with $l_1 \|P\|\|C\tilde x_o +\nu\|\leq a_6$ where $a_6$ is a bounded constant. Hence, one has $ \|{\hat{\bar\beta}}(t)\|\leq\|{\hat{\bar\beta}}(0)\|+\frac{a_6\epsilon}{\lambda}\left( 1-\exp(-\frac{\lambda}{\epsilon} t)\right)$, and it is valid to conclude that  $\hat {\bar\beta}\in \mathcal{L}_\infty$. 
   $\hfill  \blacksquare  $

  \section{Simulation Results}
  In this section, two simulation results are presented to shed some light  on  the presented theoretical discussions. In the first simulation,  a numerical example is considered and the obtained   results for MHGO-based controller are compared with the conventional HGO-based approach as well as HGO with switching gain     strategy \cite{khalilsw2009}. In the second example, simulations are carried out on a mechanical  system, and the superiority of the MHGO-based approach over conventional HGO-based   schemes and multi-observer-based approaches is  shown. 
  
  \subsection{Example 1: Underwater Vehicle }
  In this subsection, a simplified model of underwater vehicle in yaw   with dynamical equation of
  \begin{equation*}
      \ddot \psi+a \dot \psi |\dot\psi|=u
  \end{equation*}
  is selected for simulation purposes, where $\psi$ denotes the heading angle and $a$ is a positive constant. Let us assume  that only  the heading angle is   measured  and that measurement is contaminated by noise $\nu(t)$, i.e., $y=\psi+\nu(t)$. In this simulation, $a=1$ and the measurement noise, generated by Matlab uniform random number block with sampling time  $0.0001$, is in  the interval $[-0.01,0.01]$.  
  
  The control objective is to steer the heading angle to follow  the sinusoidal wave  $y_d=5+\sin(2t)$. It is obvious that the state feedback controller $u=a\dot\psi|\dot\psi|+\ddot y_d+4(\dot\psi-\dot y_d)+4(\psi-y_d)$ can force the heading angle to track $y_d$  asymptotically. Since  $\dot \psi$ is not  measurable, it should be reconstructed appropriately and fed into the above controller. As stated in the previous section,  MHGOs are capable of providing such an estimation; hence they are used in this regard.    For comparison purposes, simulations are also performed  by utilizing  the conventional HGO and  the HGO with switching   gain \cite{khalilsw2009}.  The basic idea behind  the letter observation scheme, HGO with switching gain, is to switch from a small gain to  a larger gain. More clearly, a small value for $\epsilon$ is employed in the beginning    to  get  a fast response, then it is switched to   a larger value    to avoid large steady observation errors caused  by  the measurement noise \cite{khalilsw2009}.

  The design parameters of the conventional HGO are selected as $\epsilon=0.15$, $\kappa_1=2$, $\kappa_2=1$.  Furthermore, the design parameters of the switching HGO change from $\epsilon=10^{-3}$, $\kappa_1=71$, $\kappa_2=70$ to $\epsilon=0.15$, $\kappa_1=2$, $\kappa_2=1$. Besides, it is assumed that the implemented control effort by the actuator is restricted by amplitude of $500$. For the MHGO, to be able to run the adaptive  laws  \eqref{RLS}, initial conditions of the RLS algorithm are selected as follows:  $\gamma=10^3$, {$\hat\beta_1(0)=\hat\beta_2(0)=0$}, and in turn $\hat\beta_3(0)=1-\sum\nolimits_{i=1}^{2}\hat\beta_i(0)=1$. In addition, parameters $\kappa_i$ and $\epsilon$ are set equal to the corresponding parameters of the conventional HGO.  It is clear that the  three approaches will eventually have $\epsilon=0.15$,  and in turn, all of them  will affect  the measurement noise with the same $\epsilon$, which allows us to make an accurate comparison.   Moreover, the initial conditions of the multiple observers, employed in MHGO, are required to be selected such that the initial system states, i.e., $ x(0)=[\begin{matrix} 0&0\end{matrix}]^T$ lie in their convex hull. Toward this end, three  observers are initiated from    $\hat x_1(0)=[\begin{matrix} 5&5\end{matrix}]^T,$ $\hat x_2(0)= [\begin{matrix} -5&5\end{matrix}]^T, $ and $\hat x_3(0)=[\begin{matrix} 5&-5\end{matrix}]^T$; thus $\hat x_o$ starts from  $\sum\nolimits_{i=1}^{3}\hat\beta_i(0)\hat x_i(0)=[\begin{matrix} 5&-5\end{matrix}]^T$. To make the simulation results more comparable, the initial condition of the conventional HGO and HGO with switching are set   equal to  $\hat x_o(0)=[\begin{matrix} 5&-5\end{matrix}]^T$ as well.

\begin{figure} 
\centering
        \includegraphics[ width=0.45\textwidth]{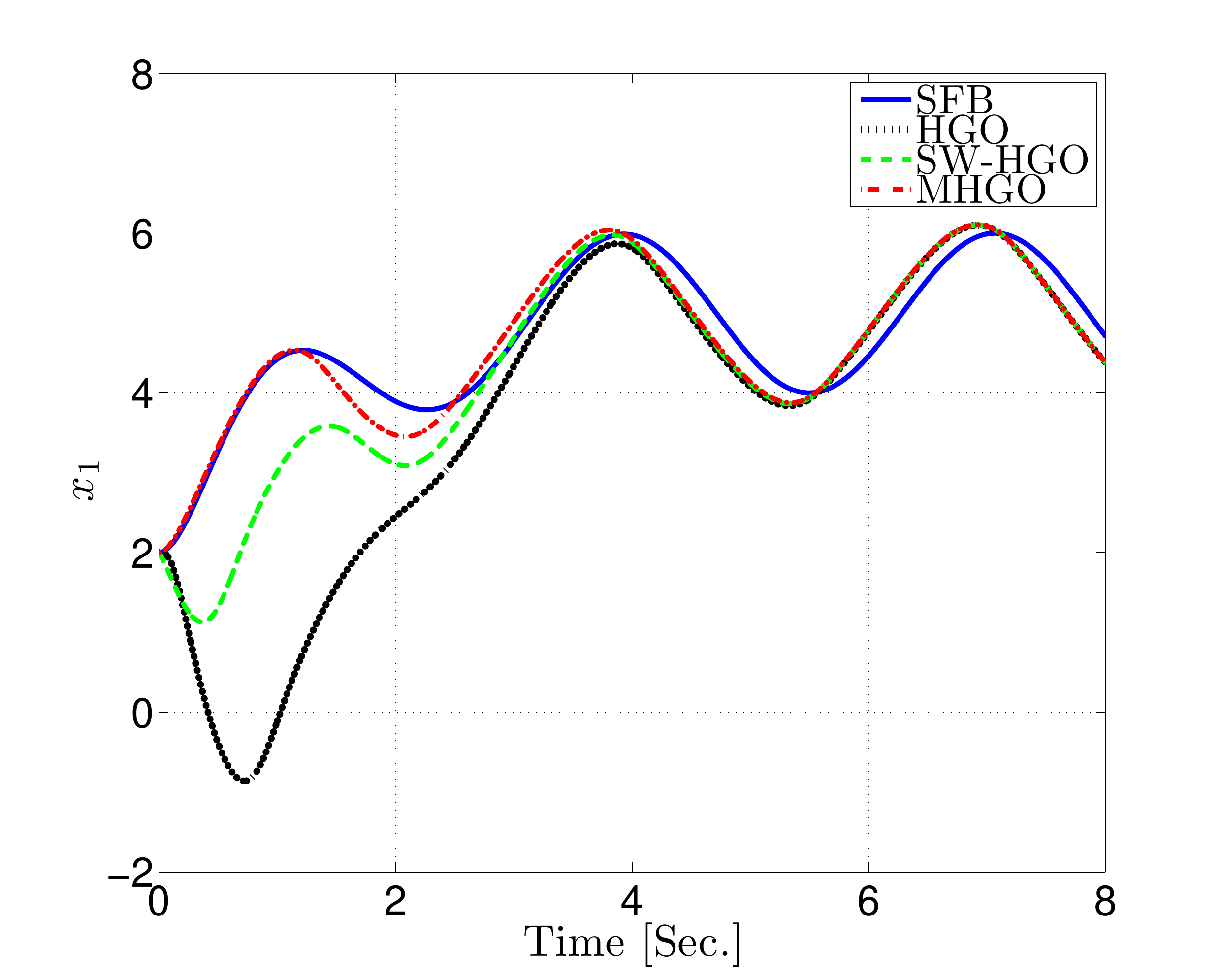}
        \caption{Evolution of $x_1(t)$ using state feedback, HGO-based, Switching HGO-based, and MHGO-based controllers.}
\label{x1}
\end{figure}

\begin{figure} 
\centering
        \includegraphics[ width=0.45\textwidth]{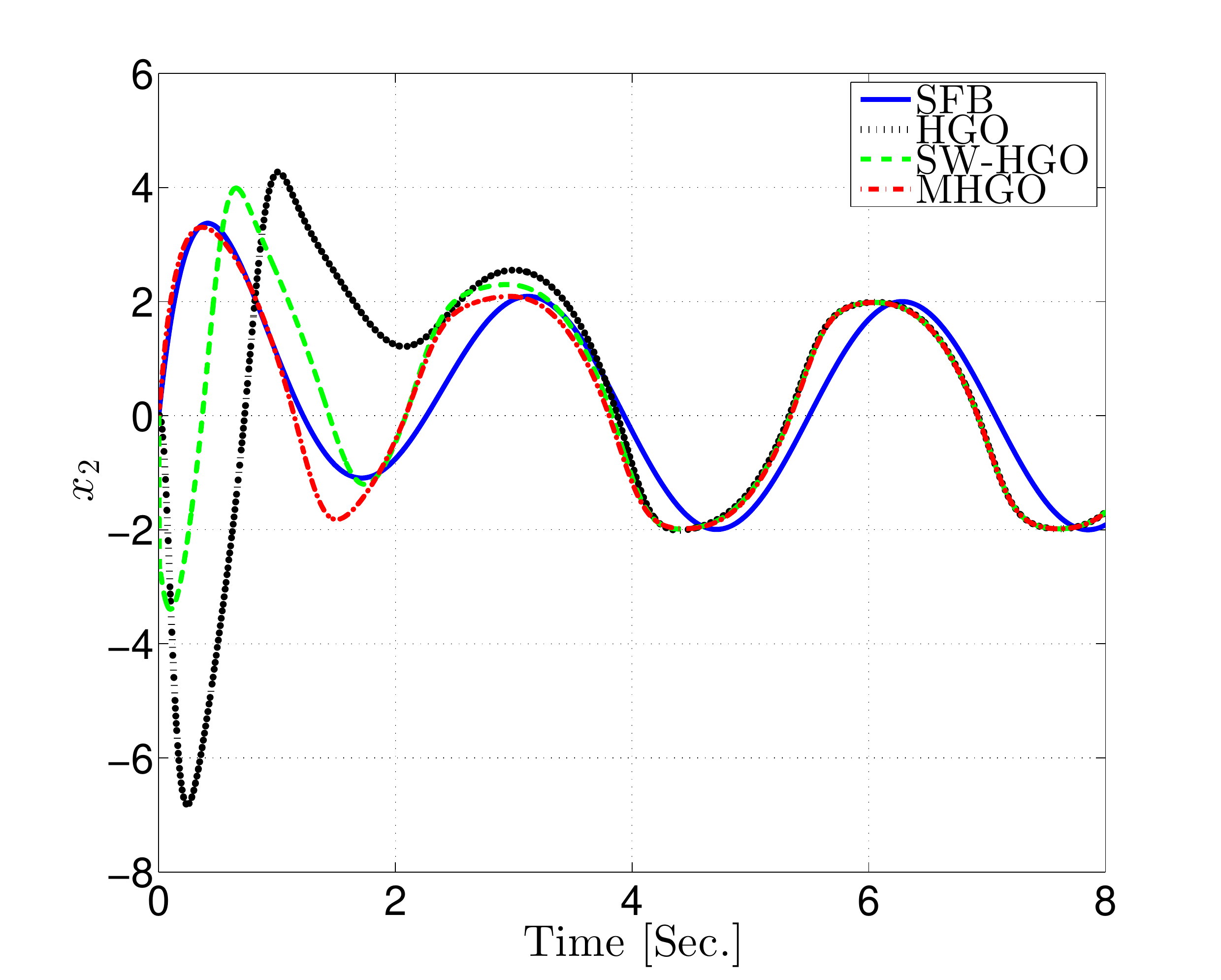}
        \caption{Evolution of $x_2(t)$ using state feedback, HGO-based, Switching HGO-based, and MHGO-based controllers.}
\label{x2}
\end{figure}

\begin{figure} 
\centering
        \includegraphics[ width=0.45\textwidth]{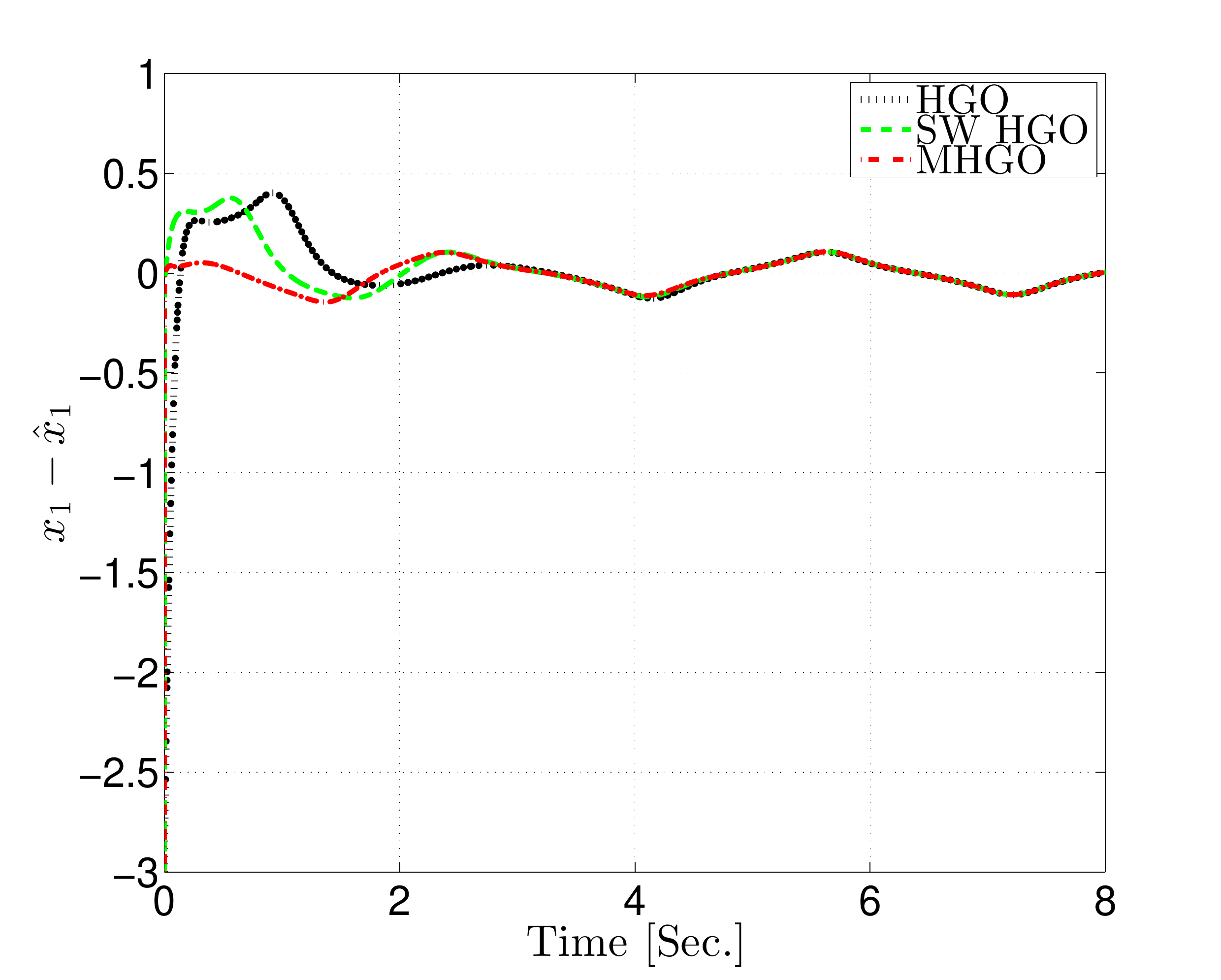}
        \caption{State estimation errors of $x_1$ using conventional  HGO,  Switching HGO and MHGO.}
\label{e1}
\end{figure}

\begin{figure} 
\centering
        \includegraphics[ width=0.45\textwidth]{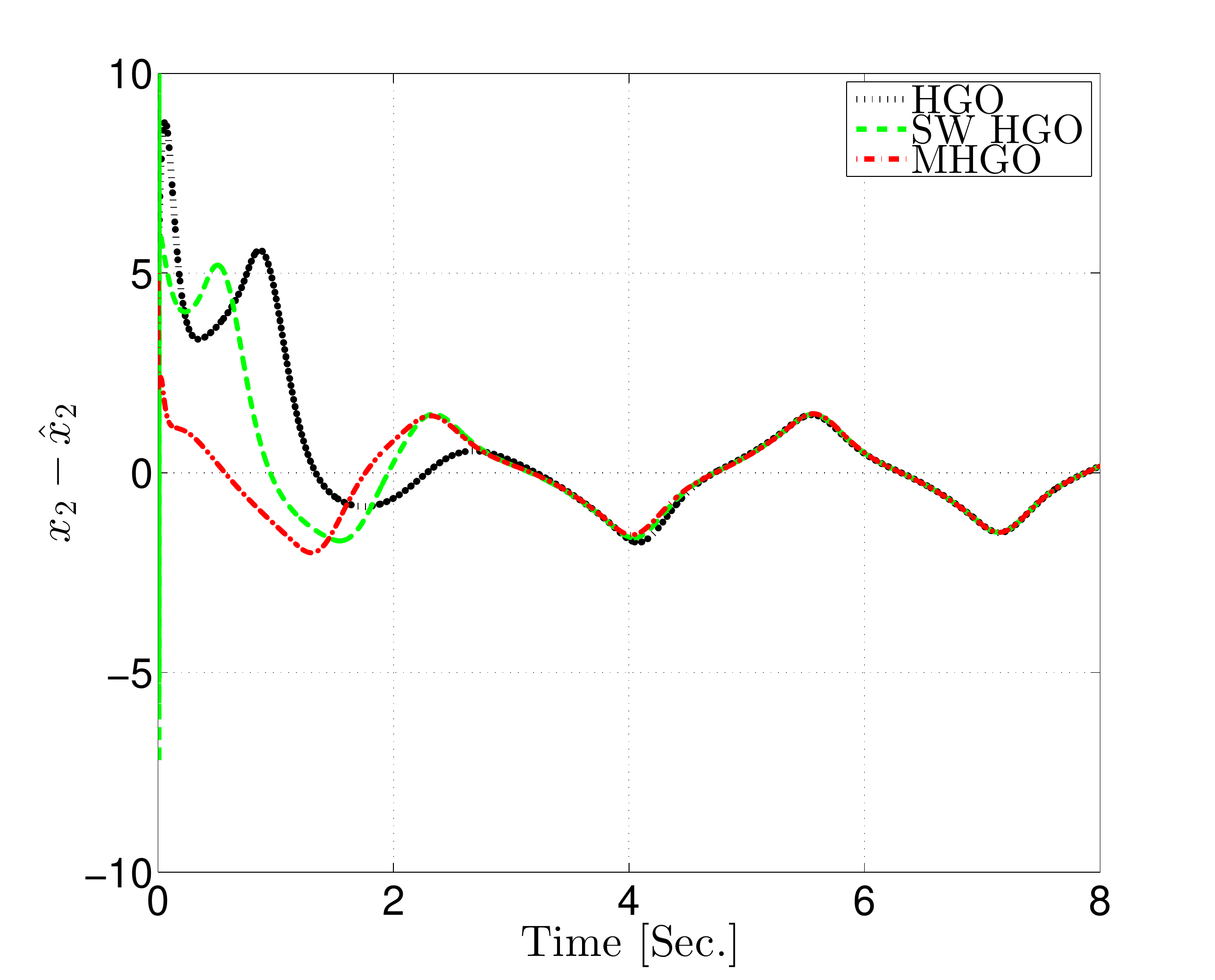}
          \caption{State estimation errors of $x_2$ using conventional  HGO,  Switching HGO and MHGO.}
\label{e2}
\end{figure}

\begin{figure} 
\centering
        \includegraphics[ width=0.45\textwidth]{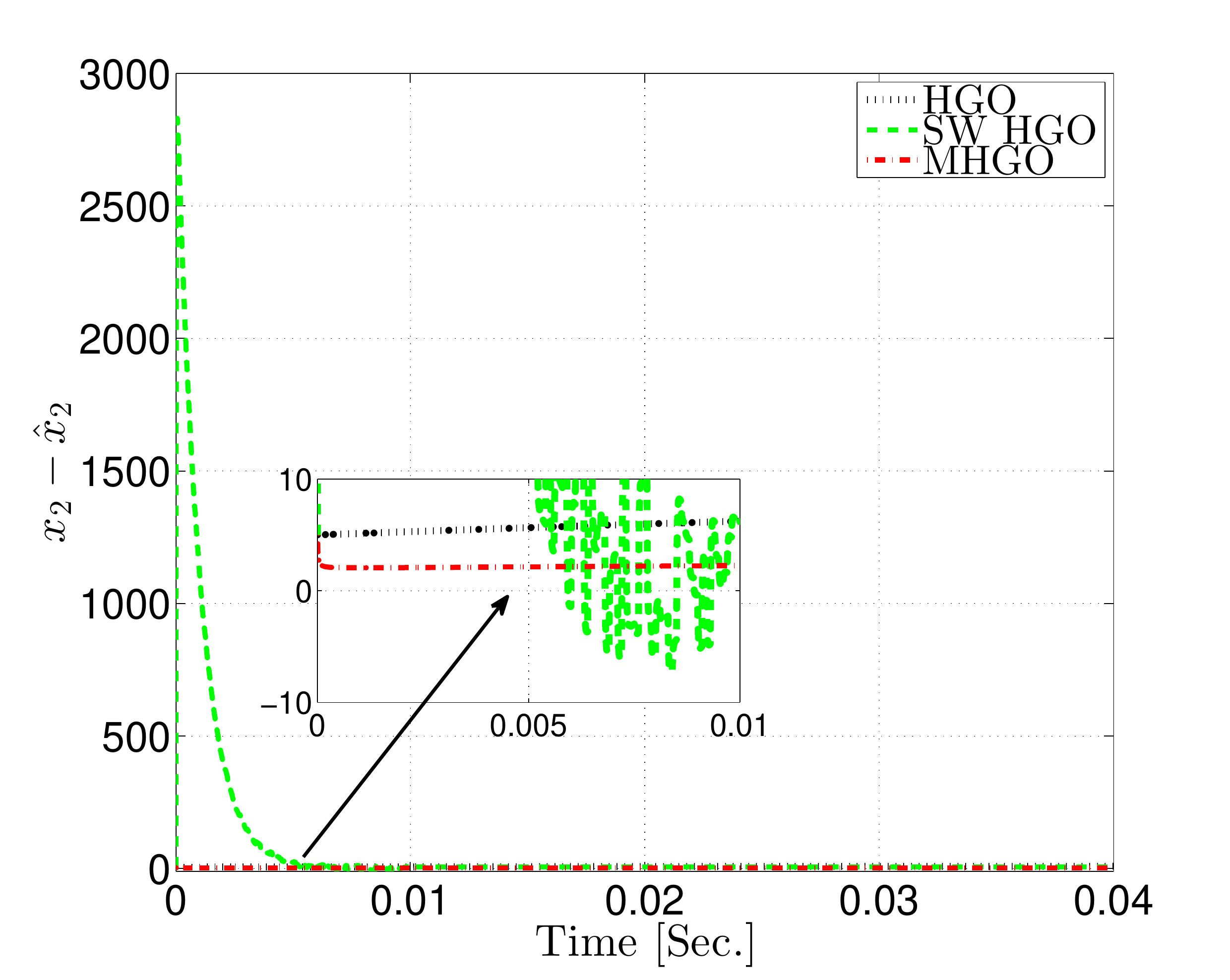}
         \caption{State estimation errors of $x_2$ during the transient phase using conventional  HGO,  Switching HGO and MHGO.}
\label{e2Z}
\end{figure}

The evolution of system states are depicted in Figs. \ref{x1} and \ref{x2}. From these two figures, it is clear that the MHGO-based controller recovers performance of the state feedback controller faster than the two others. Also,   due to the existence of the measurement noise on the system output $y(t)$, the controllers based on these observation strategies result in a bounded error. 
 The state estimation process is also presented in Fig. \ref{e1} and \ref{e2}. The obtained observation results are also in commensurate with the aforementioned discussions. These figures show that because  a small value is considered for  $\epsilon$ in the transient phase in  switching gain approach,  this scheme reconstructs system states faster than the conventional HGO. Nonetheless, as it is well-known in the control literature, such a selection yields to a large overshoot in the beginning of the estimation (see Fig. \ref{e2Z}), which is not preferable in practice. While  the obtained final estimation error for MHGO 
outperforms   the other two approaches regarding the convergence rate and the avoidance of peaks since MHGO  benefits from the advantages of using information/estimations gathered from various   sources. 
Because,  during the transient phase,  the adaptive terms $\hat\beta_i$, , presented for re-parameterizing   the observation problem and combining the estimations of each observer, assist the observation structure to result in  better estimations  (see  Remark 1).  Note that in the long run $\hat\beta_i$ converge to their final values and  since $\sum\nolimits_{i=1}^3\hat\beta_i=1$,  MHGO behaves similarly to a HGO.

  \subsection{Example 2: Connected Inverted Pendulums on Carts}
\begin{figure} 
  	\centering
  	\includegraphics[ width=0.45\textwidth]{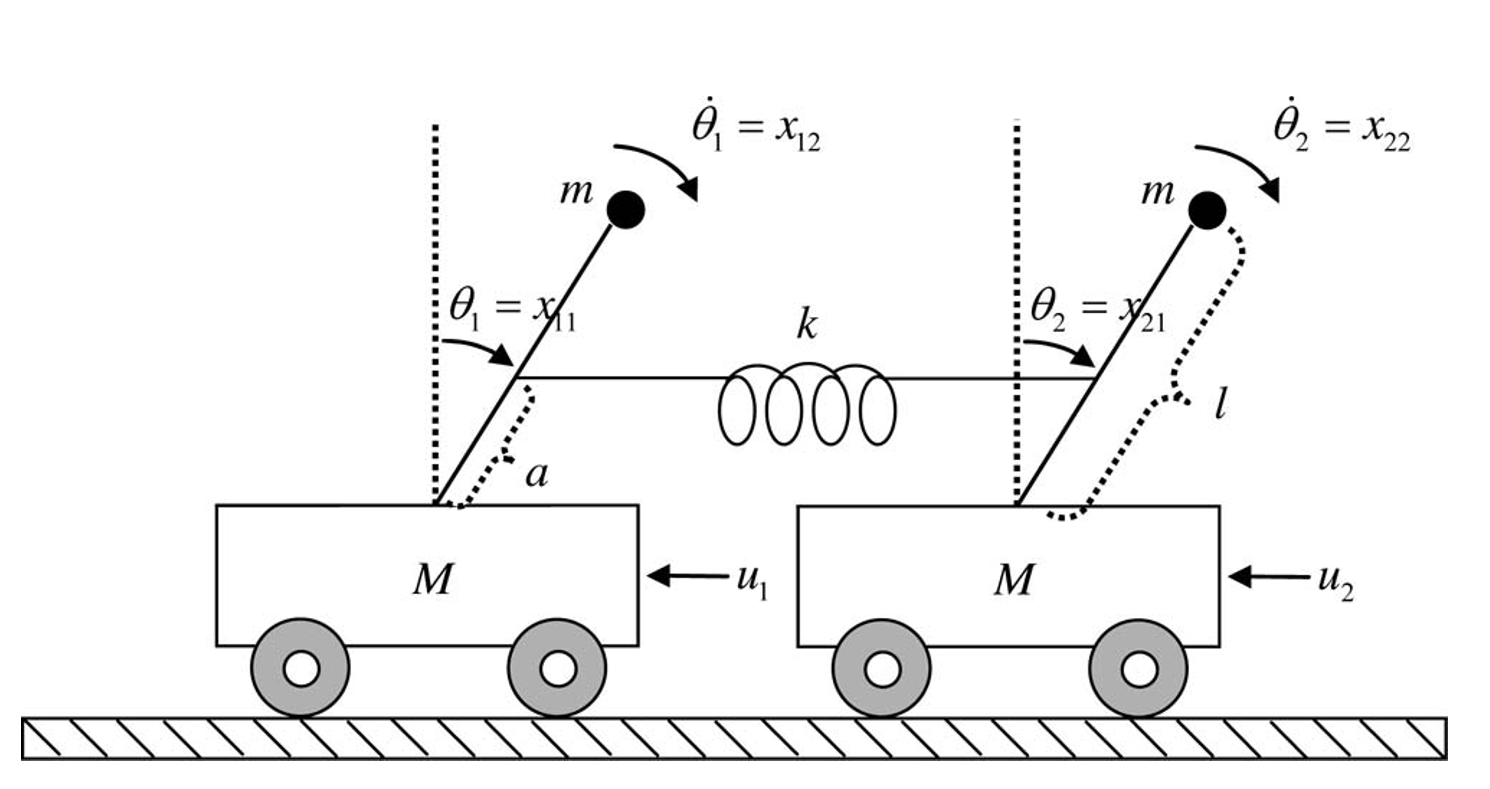}
  	\caption{A schematic of the  system considered in example 2 \cite{simulation}. }
  	\label{scheme}
  \end{figure}

  Consider two  inverted  pendulums    mounted on two carts and connected by a spring, as shown in Fig. \ref{scheme}. 
 The governing dynamical equation of this mechanical system can be expressed as \cite{simulation}, 
      \begin{equation}
      \begin{split}
          \dot x_{11}&=x_{12}\\
          \dot x_{12}&=\mathcal{F}_{11}(x)+\mathcal{F}_{12}u_1\\
            \dot x_{21}&=x_{22}\\
          \dot x_{22}&=\mathcal{F}_{21}(x)+\mathcal{F}_{22}u_2\\
          y&=[\begin{matrix}
              x_{11}&x_{21}
          \end{matrix}]^T+ \nu (t) 
               \end{split}
               \label{mech}
  \end{equation}
  where $x_{k1} $ and $x_{k2}$ ($k=1,2$) represent the vertical angle of the $k${th} pendulum  $\theta_k$  and its   angular velocity $\dot\theta_k$, respectively,      $x=[\begin{matrix}x_{11} &x_{12} &x_{21} &x_{22} \end{matrix}]^T$ is the system state vector,  $\nu(t)=[\begin{matrix}
              \nu_{1}(t)&\nu_{2}(t)
          \end{matrix}]^T $ denotes   measurement noise,   $\mathcal{F}_{k1}(x)=\left(\frac{g}{cl}-ka\frac{a-cl}{cml^2}\right)x_{k1}-\frac{m}{M}\sin(x_{k1})x^2_{k2}+ka\frac{a-cl}{cml^2}x_{j1}$ and $\mathcal{F}_{k2}(x)=\frac{1}{cml^2}$ with   $ k,j=1,2  $ and  $k\ne j$. 
                   In this simulation, the following values are considered for system parameters:  mass of pendulum $
m=1$ kg, mass of cart $M=5 $ kg, constant term $ c=\frac{m}{m+M}$, distance of the cart from the spring along with the bar $
a=0.2$ m, length of pendulum
$l=1$ m, spring constant
$
k=1$ N/m,  and gravity acceleration $
g=9.8$ m/s$^2$.

                    In order to make the closed-loop system  asymptotically stable and force its output vector to track the desired trajectory   $y_d=\left[\begin{matrix}y_{1d} & y_{2d}
          \end{matrix}\right]^T= \left[\begin{matrix}0.3\sin(t) & 0.3\cos(t) 
          \end{matrix}\right]^T$,  the state  feedback controller $$u_k=\frac{1}{\mathcal{F}_{k2}}\left(-\mathcal{F}_{k1}+\ddot y_{kd}-7(x_{k2}-\dot y_{kd})-12(x_{k1}- y_{kd})\right)$$  which is saturated outside $\left[\begin{matrix}-50 & 50
          \end{matrix}\right]$, is considered. 
          Since it was assumed that only $y_1=x_{11}+\nu_1(t) $ and $y_2=x_{21}+\nu_2(t) $ are available, system states should be reconstructed appropriately. In the sequel, the  state reconstruction process is performed using three different approaches, namely  a single HGO, multi-observer, and MHGO, and  capabilities  of  these observation strategies in recovering performance of the state feedback controller are compared. 
          
          Note that   \eqref{mech} is a multi-input multi-output (MIMO) system, hence, it is required to explain  how  the MHGO scheme can be employed for estimating the  system states. Toward this end, by using the fact that    \eqref{mech} represents a special class of MIMO systems with two subsystems in normal form, one can employ two sets of MHGO  for state estimation of the overall system. One MHGO uses $y_1$ to estimate  $x_{11}$ and $x_{12} $, and the other one employs $y_2$ for  estimating $x_{21}$ and $x_{22}$. Each set of MHGO  has $N_{k}$ HGOs  and the  parameters estimations ($\hat\beta_{k1},\hat\beta_{k2},\cdots,\hat\beta_{kN_k}$)  obtained from the  RLS algorithm \eqref{RLS}.  It is clear that   the aforementioned approach  can be easily extended to a  class of MIMO nonlinear systems consisted of more than two ($k>2$) subsystems in normal form. 

           To carry out the simulations, the  initial conditions of the system are considered as $x(0)=\begin{bmatrix} 1 & 0 & 1 & 0\end{bmatrix}^T$,  and the measurement noise vector $\nu(t)$ is generated by two    uniform random number  blocks  of Matlab Simulink with the values  restricted to the interval $[\begin{matrix}-0.02 & 0.02\end{matrix}]$ and  sampling time $0.0001$. Note that from this point forward, the  superscript  $k=1,2$ is used to denote the state estimation of  subsystem $k$. In order to investigate the performance of conventional HGO-based controller, the state variables of each subsystem are estimated using a single HGO. For that,  the design parameters and initial conditions of the $k${th} observers are  selected  as: $\kappa_{k1}=2$, $\kappa_{k2}=1$, $\epsilon_k=0.05$, $\hat x^k(0)=[\begin{matrix}3&-3\end{matrix}]^T$. In   multi-observer approach, $N_k$ HGOs are run  from various initial conditions to estimate states of the $k ${th} subsystem, and at each time instant, the performance criterion $\mu_{ki}$ {obtained from  $\dot \mu_{ki}=-\alpha\mu_{ki}+ (y-\hat y_{ki})^2:\,\,\mu_{ki}(0)=0, \alpha=0.1>0$   is checked to find the best observer} \cite{postoyan2015multi}. More clearly, the best observer for subsystem    $k$ is chosen as               $\sigma_k(t)=\text{arg}\underset{i }\min\,\,(\mu_{ki}(t))$, and in turn, the best  estimation is $\hat x^k_{mul}=\hat x_{\sigma_k}$  {(note that the observers  do not cooperate)}.
           For multi-observer approach, first we will use three observers, $N_k=3$, with initial conditions   $\hat x^k_{1}(0)=[\begin{matrix}3 &3\end{matrix}]^T$, $\hat x^k_{2}(0)=[\begin{matrix}-3& 3\end{matrix}]^T$, and $\hat x^k_{3}(0)=[\begin{matrix}3&-3\end{matrix}]^T$ to estimate the states of subsystem $k$. Moreover, the rest of design parameters of the HGOs employed in  the multi-observer method  are chosen as same as the single HGO. To be able to make a reasonable comparison between the performance recovery of the MHGO-based controller and the aforementioned methods,   the   design parameters   (i.e., $\kappa_{k1}$, $\kappa_{k2}$,  $\epsilon_k$),  the number of observers (i.e., $N_k$),  and   initial conditions (i.e., $\hat x^k_i, i=1,\cdots, N_k$) of MHGO are set equal to the  ones selected for the multi-observer scheme.  The   RLS algorithm design parameters  are considered as $P_k(0)=10^3I_{2\times2}$, $\hat\beta_{k1}=\hat\beta_{k2}=0, \hat\beta_{k3}=1$. Thus, it is clear that such a selection results in $\hat x^k_{o}(0)=\sum\nolimits_{i=1}^{3}\hat\beta_{ki}(0)\hat x^k_i(0)=\hat x^k_3(0)$. Furthermore,  $\sigma_k(0)$ in multi-observer is set as $3$; hence all the observers have the same initial condition, i.e.,  $\hat x_o^k(0)=\hat x^k_{mul}(0)=\hat x^k(0)$. 
           
\begin{figure} 
\centering 
        \includegraphics[ width=0.45\textwidth]{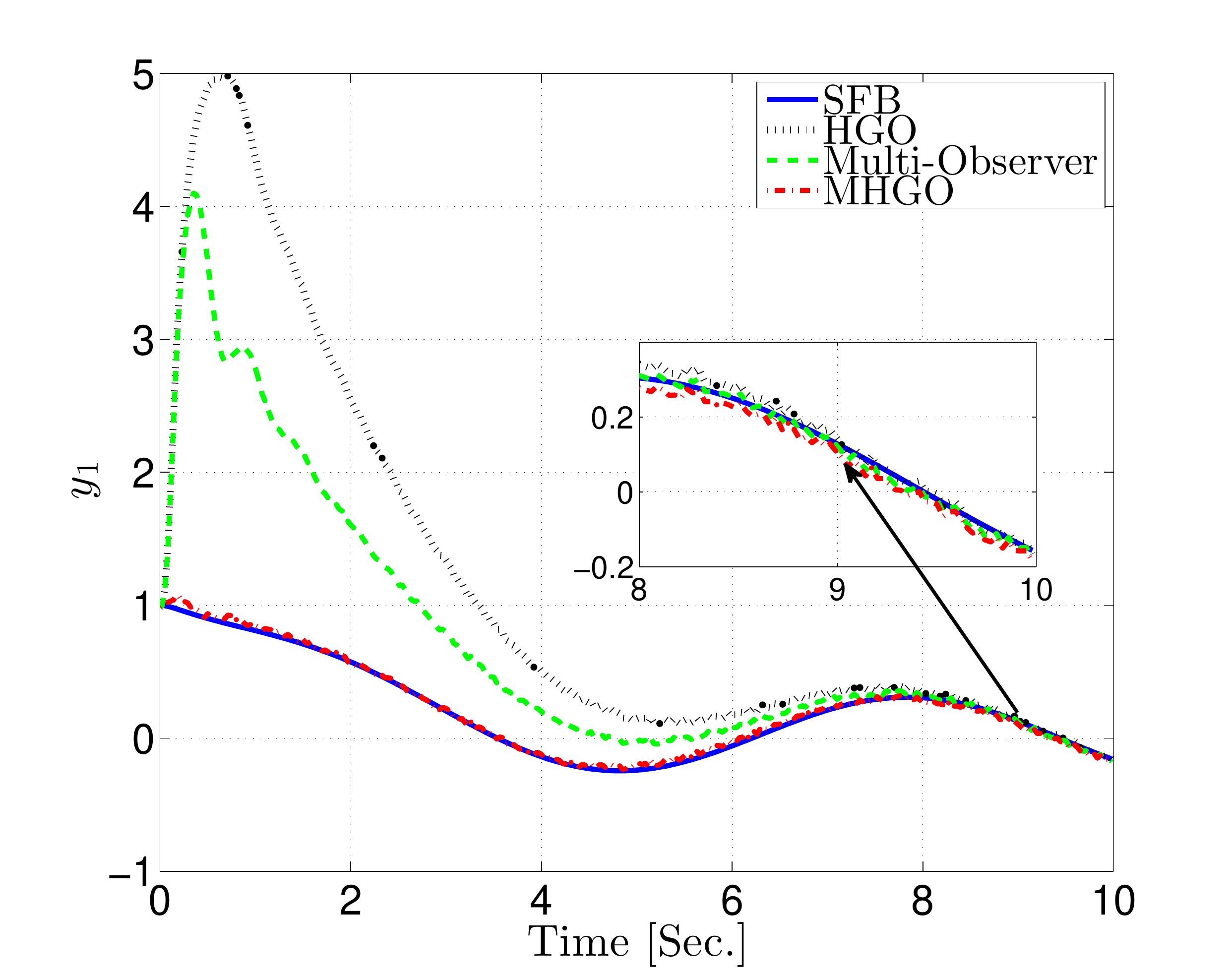}
         \caption{ Evolution  of $y_1=x_{11}+\nu_1(t)$ using  state  feedback,  HGO-based,
multi-observer-based, and MHGO-based controllers with $N_k=3$.}
\label{x11}
\end{figure}

\begin{figure} 
\centering
        \includegraphics[ width=0.45\textwidth]{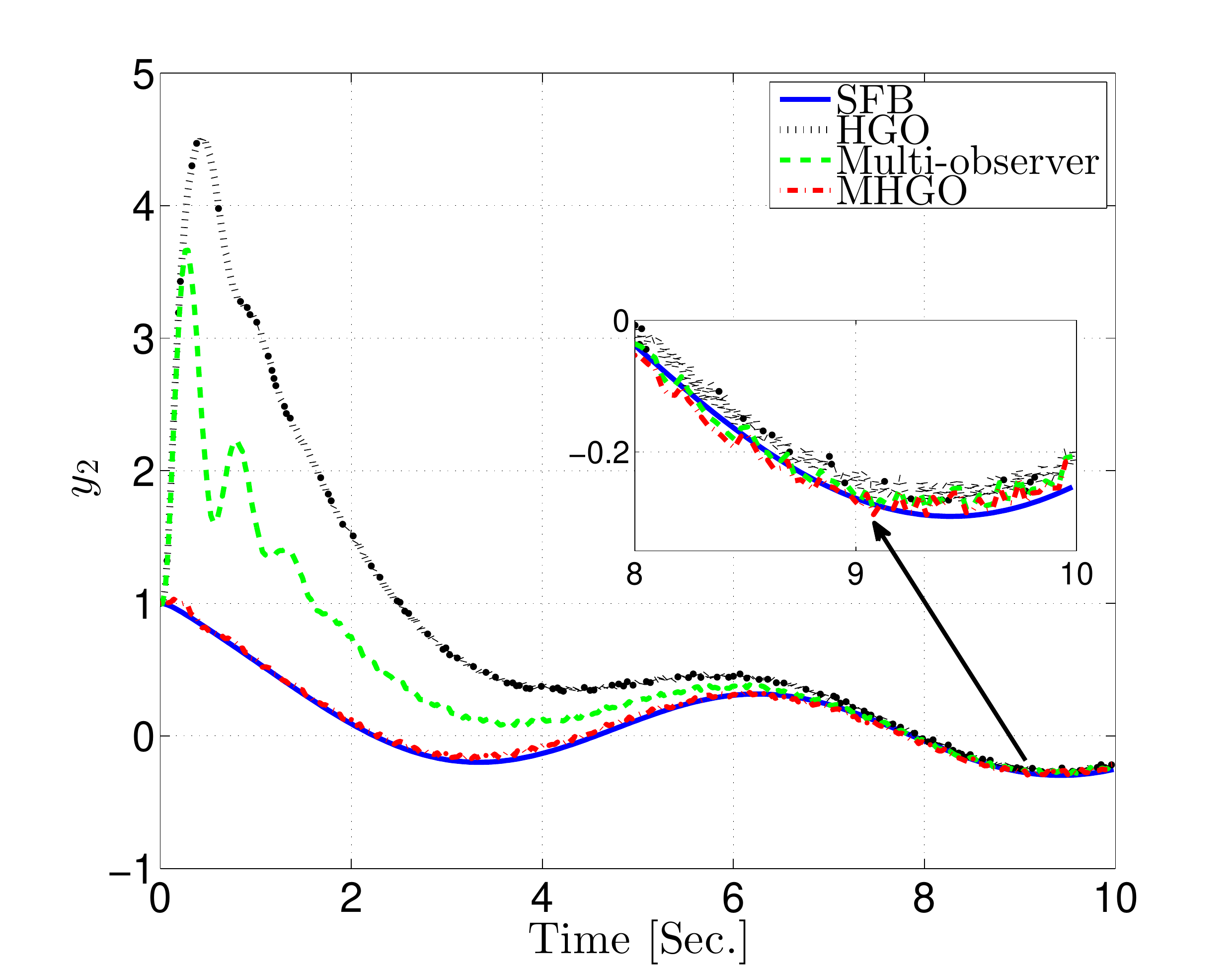}
              \caption{ Evolution  of $y_2=x_{21}+\nu_2(t)$ using  state  feedback,  HGO-based,
multi-observer-based, and MHGO-based controllers with $N_k=3$.}
\label{x21}
\end{figure}

\begin{figure} 
\centering 
        \includegraphics[ width=0.45\textwidth]{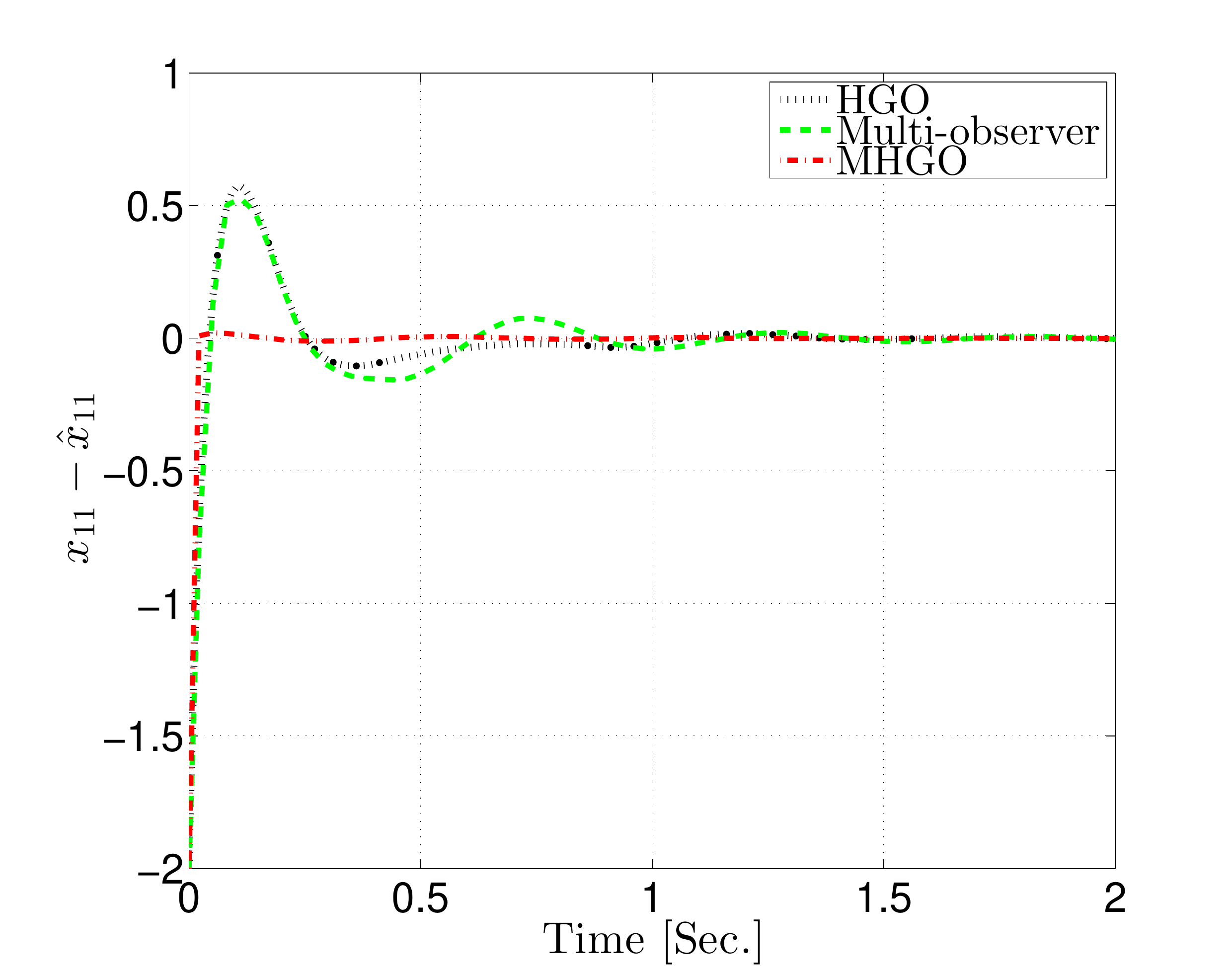}
 \caption{Estimation error of first state, $x_{11}$, using conventional  HGO,  multi-observer,  and MHGO with $N_k=3$.}
 \label{ex11}
\end{figure}

\begin{figure} 
\centering 
        \includegraphics[ width=0.45\textwidth]{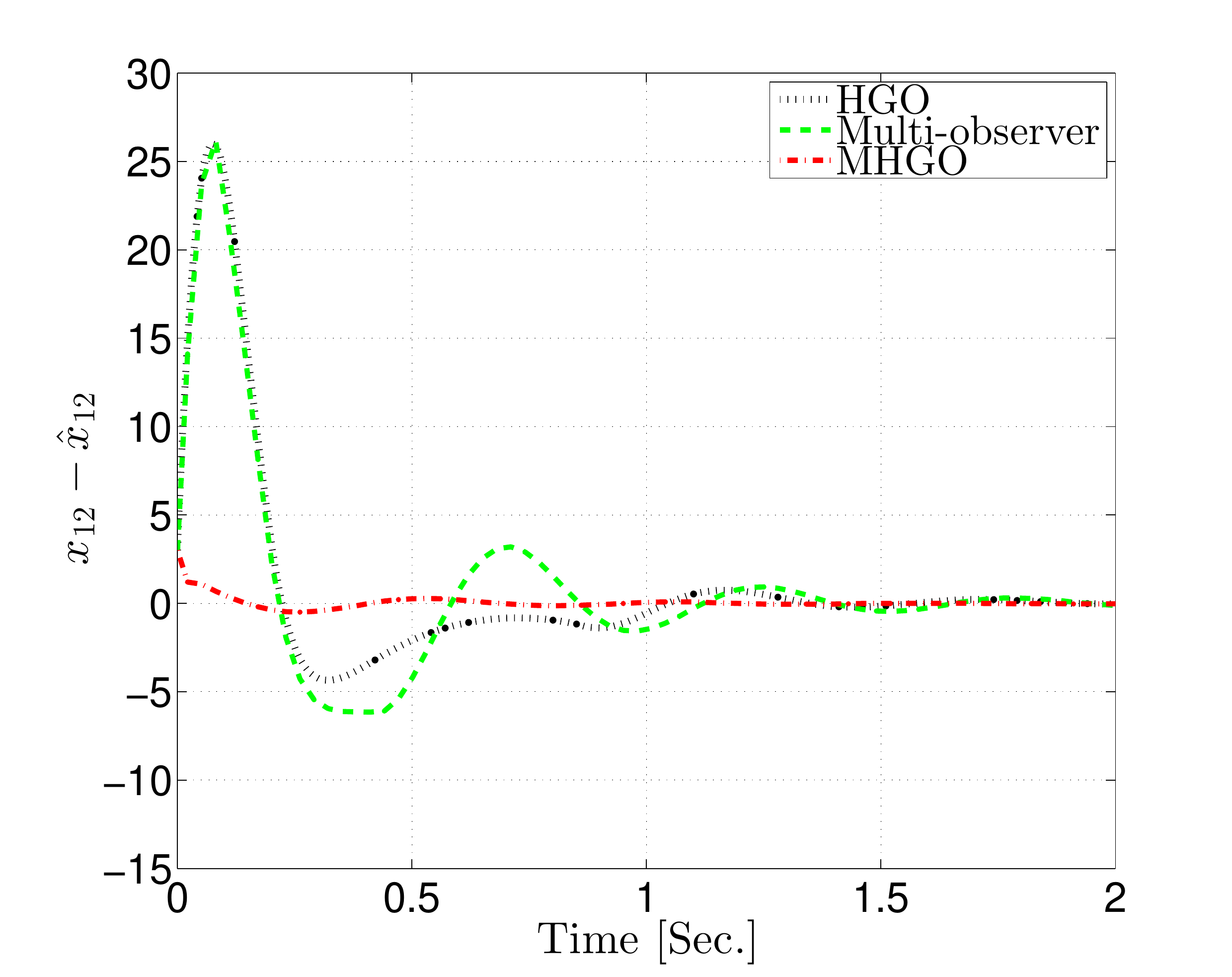}
         \caption{Estimation error of second state, $x_{12}$, using conventional  HGO,  multi-observer, and MHGO with $N_k=3$.}
\label{ex12}
\end{figure}

           The  performance of the system states under the state feedback controller as well as the discussed output feedback controllers are illustrated in   Figs. \ref{x11} and \ref{x21}.  These  two figures obviously show that although the multi-observer-based controller yields a better response in comparison to the conventional HGO-based method, it does not outperform   the MHGO-based controller. In other words, the MHGO-based controller  has recovered the performance of the      state feedback controller much faster than the other two  methodologies. As discussed earlier and it is obvious  from the zoomed parts in F-gs. \ref{x11} and \ref{x21},   performance of the state feedback controller cannot be recovered perfectly due to the existence of measurement noise $\nu(t)$.  Because of the space limitation and the fact  the two subsystems behave  similarly, only the  observation errors of the first subsystem  are   provided (Figs. \ref{ex11} and \ref{ex12}). According to these figures,  one can easily see that the state estimation errors of MHGO strategy converge to a small neighborhood of the origin rapidly.

           

\begin{figure} 
\centering
        \includegraphics[ width=0.45\textwidth]{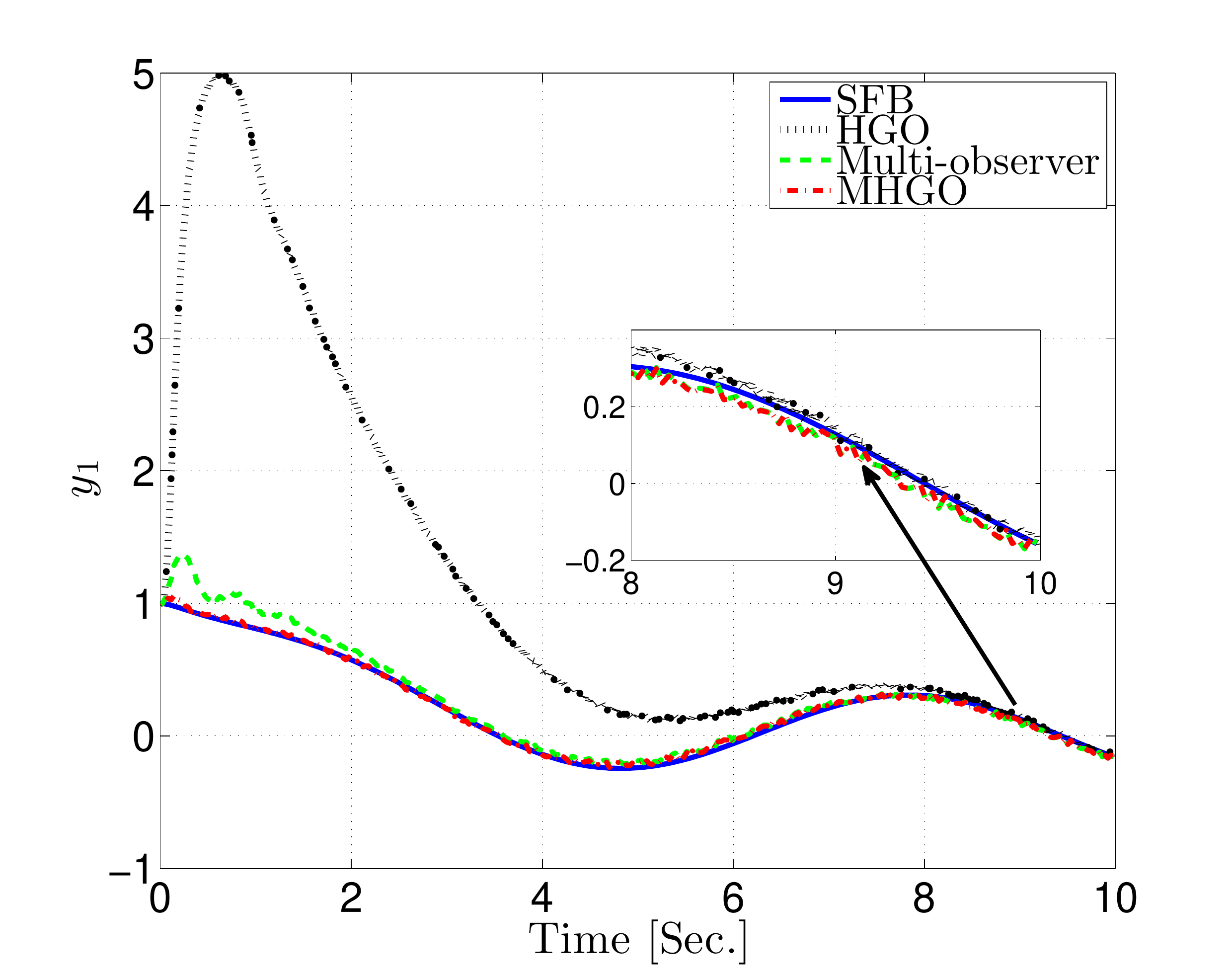}
           \caption{ Evolution  of $y_1=x_{11}+\nu_1(t)$ using  state  feedback,  HGO-based,
multi-observer-based, and MHGO-based controllers with $N_k=81$.}
\label{nx11}
\end{figure}

\begin{figure} 
\centering
        \includegraphics[ width=0.45\textwidth]{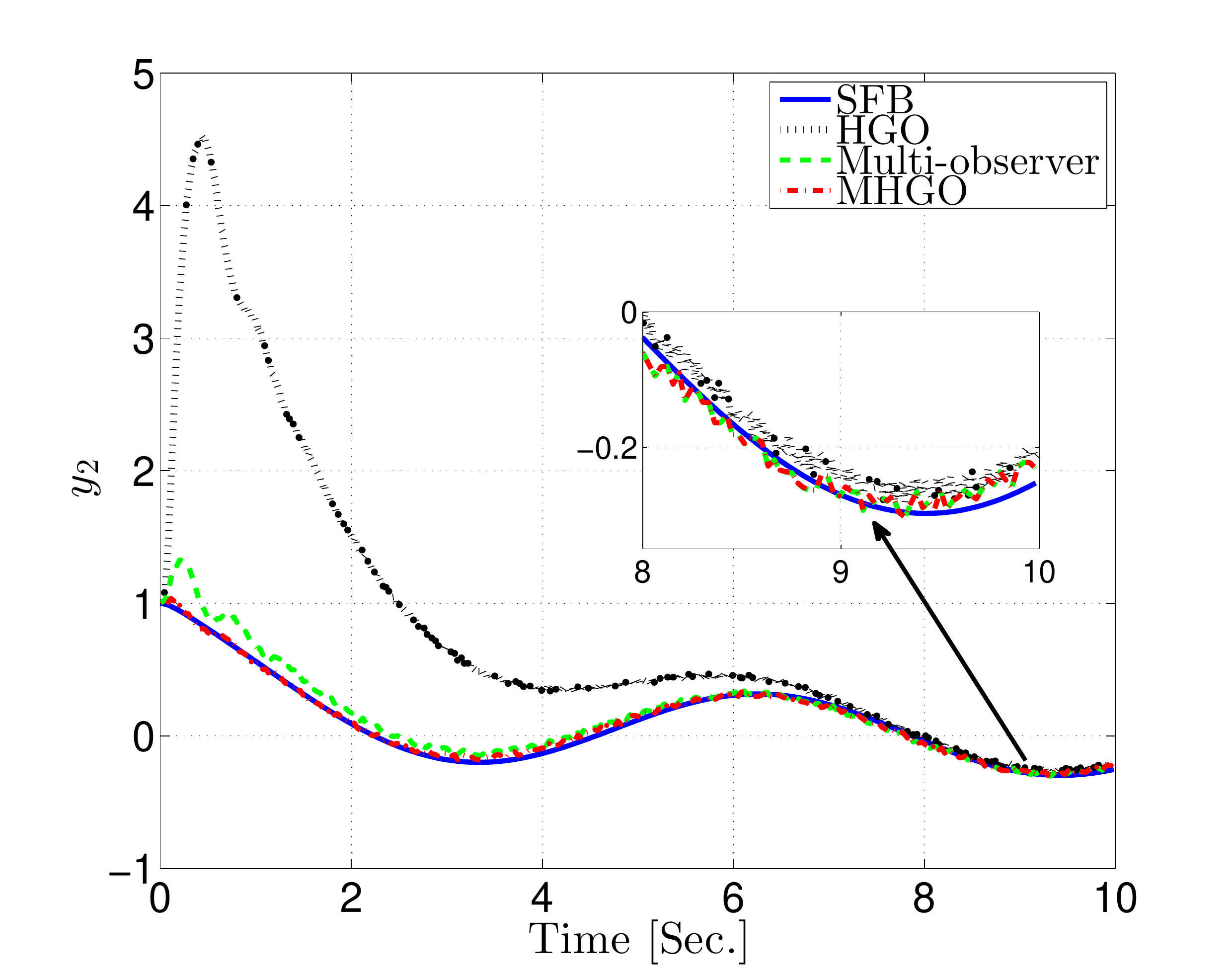}
          \caption{ Evolution  of $y_2=x_{21}+\nu_2(t)$ using  state  feedback,  HGO-based,
multi-observer-based, and MHGO-based controllers with $N_k=81$.}
\label{nx21}
\end{figure}

\begin{figure} 
\centering 
        \includegraphics[ width=0.45\textwidth]{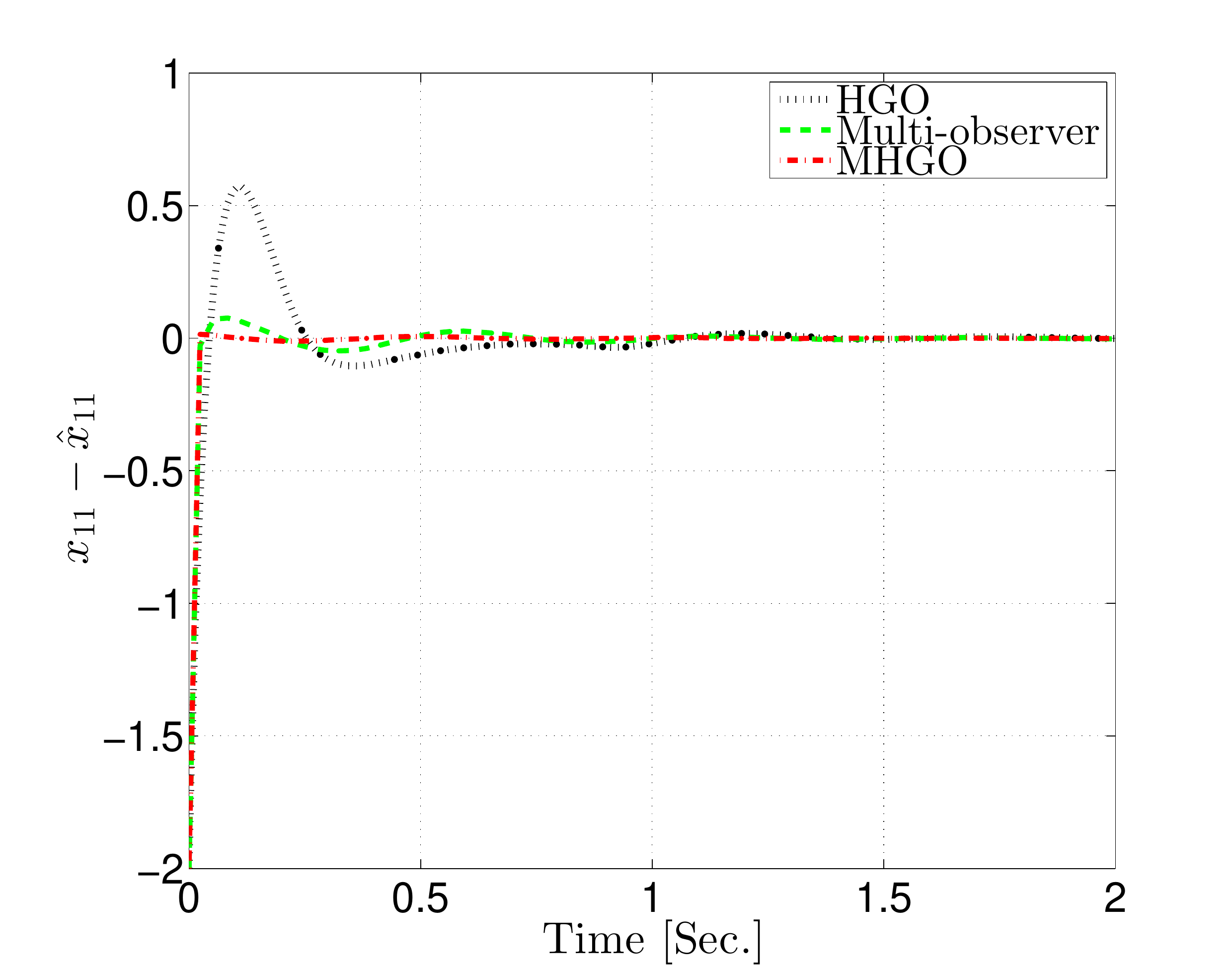}
        \caption{Estimation error of first state, $x_{11}$, using conventional  HGO,  multi-observer,  and MHGO with $N_k=81$.}
\label{nex11}
\end{figure}

\begin{figure} 
\centering 
        \includegraphics[ width=0.45\textwidth]{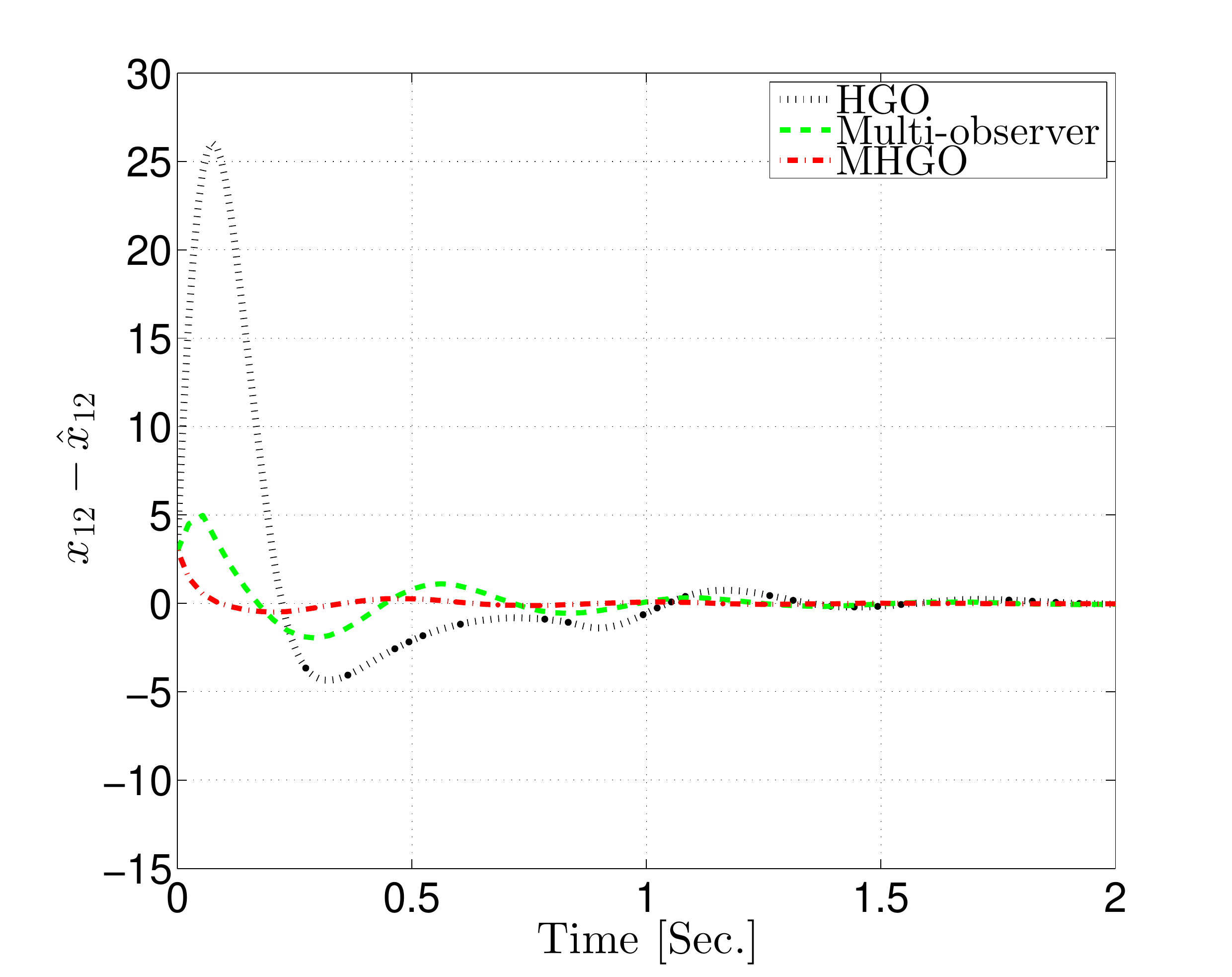}
         \caption{Estimation error of second  state, $x_{21}$, using conventional  HGO,  multi-observer,  and MHGO with $N_k=81$.}
\label{nex12}
\end{figure}

            {As mentioned in  }  \cite{postoyan2015multi}, \cite{mulobs},  {  the multi-observer approach provides better results  when the number of observers is increased. This implies that to assure that at least one of the models is sufficiently close to the plant, a  quite large number of models is required.  }  Hence,   to make a more comprehensive comparison and show that the MHGO-based controller can outperform   the multi-observer-based controller even when the number of observers are increased, another simulation with $N_k=81$ is carried out. In this simulation, all the design parameters are chosen the same as the previous scenario.
            For providing the initial conditions of observers for  subsystem $k$, four points  $[\begin{matrix}3&3\end{matrix}]^T    $, $[\begin{matrix}3&-3\end{matrix}]^T    $, $[\begin{matrix}-3&3\end{matrix}]^T  $, and $[\begin{matrix}-3&-3\end{matrix}]^T    $   are considered as the vertices of  uncertainty region $\mathcal{K}_k$, within which the initial conditions of subsystem $l$ lies. Then,   $\mathcal{K}_k$ is sampled uniformly to obtain $81$ initial conditions for each set of observers. Simulation results are presented in Figs. \ref{nx11}-\ref{nex12}.  Figs. \ref{nx11} and \ref{nx21} clearly show  that  performance of the multi-observer-based approach is improved in comparison to the multi-observer with $N_k=3$ (see Figs. \ref{x11} and \ref{x21}); however, it is computationally more expensive than that case.
          These  figures also demonstrate that the MHGO-based control approach results in a better performance in this scenario as well. It is worth noting that in these simulations the MHGO-based controller with $N_k=3$ also outperforms    the multi-observer-based controller with $N_k=81$, which demonstrates the superiority of the MHGO-based strategy (refer to  Figs. \ref{x11}, \ref{x21}, \ref{nx11},   \ref{nx21}). The obtained observations errors are also shown in Figs. \ref{nex11}  and \ref{nex12}, based on which it is clear that the MHGO scheme forces the observation errors to tend to a small neighborhood of the origin faster than other methods. 
            
            

 \section{Conclusions}
This paper investigates the analysis of the performance of MHGO-based controllers when the system output is contaminated by measurement noise. It is well-known that in this case, increasing the  gain of a single observer    deteriorates the steady state bound of estimation error. Hence, one cannot choose an arbitrarily large gain    to speed up the transient response, which is necessary for control purposes.  The  MHGO     utilizes  the observations obtained from various sources and introduces new design parameters. In turn, it provides a suitable tool for solving the  aforementioned  trade-off   in conventional HGO. The necessary conditions under which  such a  structure was capable of recovering performance of the state feedback controllers in the presence of measurement noise were  derived, and the stability of the closed-loop system was investigated. Two simulations were carried out for comparison purposes and  to validate the theoretical discussions. 

Several simulations performed on various dynamical systems have shown that the MHGO-based control strategy has a larger region of attraction in comparison to the conventional  HGO. The focus of our future work is to   investigate this problem further  and provide a rigorous proof for this claim. 


 

\begin{thebibliography}{100} 

\bibitem{RobustnessObs}
J.~C. Doyle and G.~Stein, ``Robustness with observers,'' \emph{IEEE Trans.
  Autom. Contr.}, vol.~24, no.~4, pp. 607--611, 1979.

\bibitem{Kalman}
R.~E. Kalman, ``A new approach to linear ﬁltering and prediction problems,''
  \emph{Fluids Eng.}, vol.~82, no.~1, p. 35–45, 1960.

\bibitem{Cybern2018}
Z.~Xing, Y.~Xia, L.~Yan, K.~Lu, and Q.~Gong, ``Multisensor distributed weighted
  kalman filter fusion with network delays, stochastic uncertainties,
  autocorrelated, and cross-correlated noises,'' \emph{IEEE Trans. Syst. Man.
  Cybern.: Syst.}, vol.~48, no.~5, p. 716–726, 2018.

\bibitem{Cybern2017}
B.~Chen, G.~Hu, D.~W.~C. Ho, W.~A. Zhang, and L.~Yu, ``Distributed robust
  fusion estimation with application to state monitoring systems,'' \emph{IEEE
  Trans. Syst. Man. Cybern.: Syst.}, vol.~47, no.~11, pp. 2994--3005, 2017.

\bibitem{obs1}
S.~Battilotti, ``Robust observer design under measurement noise with gain
  adaptation and saturated estimate,'' \emph{Automatica}, vol.~81, pp. 75--86,
  2017.

\bibitem{Khalil}
H.~K. Khalil and J.~Grizzle, \emph{Nonlinear systems}.\hskip 1em plus 0.5em
  minus 0.4em\relax Prentice hall New Jersey, 1996.

\bibitem{atassi1999separation}
A.~N. Atassi and H.~K. Khalil, ``A separation principle for the stabilization
  of a class of nonlinear systems,'' \emph{IEEE Trans. Autom. Control},
  vol.~44, no.~9, pp. 1672--1687, 1999.

\bibitem{cyb}
W.~He, A.~O. David, Z.~Yin, and C.~Sun, ``Neural network control of a robotic
  manipulator with input deadzone and output constraint,'' \emph{IEEE Trans.
  Syst., Man, Cybern.: Syst.}, vol.~46, no.~6, pp. 759--770, 2016.

\bibitem{4}
K.~Esfandiari, F.~Abdollahi, and H.~A. Talebi, ``Adaptive control of uncertain
  nonaffine nonlinear systems with input saturation using neural networks,''
  \emph{IEEE Trans. Neural Netw. Learn. Syst.}, vol.~26, no.~10, pp.
  2311--2322, 2015.

\bibitem{Mehran}
M.~Shakarami, K.~Esfandiari, M.~A. Shamsi, and M.~B. Menhaj, ``High-gain
  observer-based identification scheme for estimation of physical parameters of
  synchronous generators,'' \emph{24th Iranian Conference on Electrical
  Engineering (ICEE)}, pp. 1422--1427, 2016.

\bibitem{esfandiari1992output}
F.~Esfandiari and H.~K. Khalil, ``Output feedback stabilization of fully
  linearizable systems,'' \emph{Int. J. Control}, vol.~56, no.~5, pp.
  1007--1037, 1992.

\bibitem{fuzzyhgo1}
C.~Li, S.~Tong, and W.~Wang, ``Fuzzy adaptive high-gain-based observer
  backstepping control for siso nonlinear systems,'' \emph{Info. Sci.}, vol.
  181, pp. 2405--2421, 2011.

\bibitem{fuzzyhgo2}
C.~Ren, S.~Tong, and Y.~Li, ``Fuzzy adaptive high-gain-based observer
  backstepping control for siso nonlinear systems with dynamical
  uncertainties,'' \emph{Nonlinear Dyn.}, vol.~67, no.~2, pp. 941--955, 2012.

\bibitem{ICEE2015}
K.~Esfandiari, F.~Abdollahi, and H.~A. Talebi, ``Adaptive output feedback
  tracking control for nonaffine nonlinear systems,'' in \emph{2015 23rd
  Iranian Conference on Electrical Engineering}.\hskip 1em plus 0.5em minus
  0.4em\relax IEEE, 2015, pp. 976--981.

\bibitem{Narendra}
K.~S. Narendra and A.~Annaswamy, \emph{Stable adaptive systems}.\hskip 1em plus
  0.5em minus 0.4em\relax Prentice hall New Jersey, 1989.

\bibitem{iet}
K.~Esfandiari, F.~Abdollahi, and H.~A. Talebi, ``Stable adaptive output
  feedback controller for a class of uncertain non-linear systems,'' \emph{IET
  Control Theory, Appl.}, vol.~9, no.~9, pp. 1329--1337, 2015.

\bibitem{NN}
K.~Esfandiari, F.~Abdollahi, and H.~A. Talebi, ``Adaptive near-optimal neuro controller for continuous-time nonaffine
  nonlinear systems with constrained input,'' \emph{Neural Netw.}, vol.~93, pp.
  195--204, 2017.

\bibitem{NaderCDC}
K.~S. Narendra and K.~Esfandiari, ``Adaptive control of linear periodic systems
  using multiple models,'' \emph{2018 IEEE Conference on Decision and Control
  (CDC)}, pp. 589--594, 2018.

\bibitem{Switch}
K.~S. Narendra and O.~A. Driollet, ``Adaptive control using multiple models,
  switching, and tuning,'' in \emph{Adaptive Systems for Signal Processing,
  Communications, and Control Symposium}.\hskip 1em plus 0.5em minus
  0.4em\relax IEEE, 2000, pp. 159--164.

\bibitem{postoyan2015multi}
R.~Postoyan, M.~H. Hamid, and J.~Daafouz, ``A multi-observer approach for the
  state estimation of nonlinear systems,'' in \emph{2015 54th IEEE Conference
  on Decision and Control (CDC)}.\hskip 1em plus 0.5em minus 0.4em\relax IEEE,
  2015, pp. 1793--1798.

\bibitem{mulobs}
M.~S. Chong, D.~Nešić, R.~Postoyan, and L.~Kuhlmann, ``Parameter and state
  estimation of nonlinear systems using a multi-observer under the supervisory
  framework,'' \emph{IEEE Trans. Autom. Control}, vol.~60, no.~9, pp.
  2336--2349, 2015.

\bibitem{CDC}
M.~Shakarami, K.~Esfandiari, A.~A. Suratgar, and H.~A. Talebi, ``On the peaking
  attenuation and transient response improvement of high-gain observers,''
  \emph{2018 IEEE Conference on Decision and Control (CDC)}, pp. 577--582,
  2018.

\bibitem{second}
Z.~Han and K.~S. Narendra, ``New concepts in adaptive control using multiple
  models,'' \emph{IEEE Trans. Autom. Control}, vol.~57, no.~1, p. 78–89,
  2012.

\bibitem{khalil1}
L.~K. Vasiljevica and H.~K. Khalil, ``Error bounds in differentiation of noisy
  signals by high-gain observers,'' \emph{Control Syst. Letters}, vol.~57,
  no.~10, pp. 856--862, 2001.

\bibitem{tradeoff}
H.~Kwakernaak and R.~Sivan, \emph{Linear Optimal Control Systems}.\hskip 1em
  plus 0.5em minus 0.4em\relax New York: Wiley-Interscience, 1972.

\bibitem{khalilsw2009}
J.~H. Ahrens and H.~K. Khalil, ``High-gain observers in the presence of
  measurement noise: A switched-gain approach,'' \emph{Control, Syst.,
  Letters}, vol.~45, no.~4, pp. 936--943, 2009.

\bibitem{bakelman2012convex}
I.~J. Bakelman, \emph{Convex analysis and nonlinear geometric elliptic
  equations}.\hskip 1em plus 0.5em minus 0.4em\relax Springer Science \&
  Business Media, 2012.

\bibitem{simulation}
W.~Y. Wang, Y.~H. Chien, and T.~T. Lee, ``Observer-based t–s fuzzy control
  for a class of general nonaffine nonlinears systems using generalized
  projection-update laws,'' \emph{IEEE Trans. Fuzzy Syst.}, vol.~19, no.~3, pp.
  493--504, 2011.

\end{thebibliography}
\end{document}